\documentclass[11pt,a4paper,notitlepage]{article}

\usepackage{amsmath}
\usepackage{amsfonts}
\usepackage{hyperref}
\usepackage{setspace}
\usepackage{fullpage}
\usepackage{braket}
\usepackage{verbatim}
\usepackage{lscape}
\usepackage{tabularx}
\usepackage{tikz-cd}

\usepackage{upgreek}
\usepackage{pgfplots}
\usepackage{enumitem}
\usepackage{amsthm,thmtools}
\usepackage{extarrows}
\usepackage[colorinlistoftodos]{todonotes}
\usepackage{xypic}
\usepackage{graphicx}
\usepackage{graphics,amssymb}
\usepackage{dsfont}
\usepackage{mathtools}
\usepackage{url}
\usepackage{fancybox}
\usepackage{dsfont}
\usepackage{bbding}
\usepackage{titlesec}
\usepackage[toc,page]{appendix}
\usepackage[titles]{tocloft}
\usepackage{pigpen}
\usepackage{setspace}
\usepackage{hyperref}
\usepackage{bbm}
\usepackage{authblk}
\usepackage{tikz}
\usepackage{caption}
\usetikzlibrary{decorations.pathreplacing,calligraphy}
\usepackage{ragged2e}
\usepackage{wrapfig}

\definecolor{palesilver}{rgb}{0.79, 0.75, 0.73}
\definecolor{spanishgray}{HTML}{A39595}
\definecolor{darkslategray}{HTML}{2F4F4F}

\theoremstyle{plain}

\newtheorem{theorem}{Theorem}[section]
\newtheorem{Theorem}{Theorem}
\newtheorem{lemma}[theorem]{Lemma}
\newtheorem{proposition}[theorem]{Proposition}

\newtheorem*{theorem*}{Theorem}
\newtheorem*{corollary*}{Corollary}
\newtheorem*{proposition*}{Proposition}
\newtheorem{Proposition}{Proposition}

\theoremstyle{definition}
\newtheorem{remark}[theorem]{Remark}

\newcommand\radiusb{0.08}
\newcommand\radiusc{0.2}

\newcommand\anglezerozero{360*(0+0*3)/9}
\newcommand\anglezeroone{360*(0+1*3)/9)}
\newcommand\anglezerotwo{360*(0+2*3)/9)}

\newcommand\angleonezero{360*(1+0*3)/9}
\newcommand\angleoneone{360*(1+1*3)/9)}
\newcommand\angleonetwo{360*(1+2*3)/9)}

\newcommand\angletwozero{360*(2+0*3)/9}
\newcommand\angletwoone{360*(2+1*3)/9)}
\newcommand\angletwotwo{360*(2+2*3)/9)}

\newcommand{\Q}{\mathbb{Q}}
\newcommand{\R}{\mathbb{R}}
\newcommand{\Z}{\mathbb{Z}}

\newcommand{\SL}{\mathrm{SL}}

\newcommand{\Addresses}{{\bigskip \footnotesize
\textsc{Lee}
 
 \textsc{Mathematics Institute, University of Warwick,} 
 
 \textsc{Coventry CV4 7AL, UK} 
 
 \textit{E-mail address}: \texttt{jungwon.lee@warwick.ac.uk}

\medskip

\textsc{Palvannan}
 
\textsc{Department of Mathematics, Indian Institute of Science,  CV Raman Road,}

\textsc{Bangalore 560012, India} 

 \textit{E-mail address}: \texttt{bharathwaj@iisc.ac.in}
}}

\date{}
\titleformat{\section}{\bfseries\large}{\S\thesection.  }{0.1pt}{}
\titleformat{\subsection}{\bfseries}{\S \thesubsection.  }{0.1pt}{}
\titleformat{\subsubsection}{\bfseries}{\S \thesubsubsection.  }{0.1pt}{}
\title{An ergodic approach towards an equidistribution result of Ferrero--Washington}
\author{Jungwon Lee and Bharathwaj Palvannan}

\makeatletter
\newcommand{\subjclass}[2][2020]{%
  \let\@oldtitle\@title%
  \gdef\@title{\@oldtitle\footnotetext{#1 \emph{Mathematics subject classification.} #2}}%
}
\newcommand{\keywords}[1]{%
  \let\@@oldtitle\@title%
  \gdef\@title{\@@oldtitle\footnotetext{\emph{Key words and phrases.} #1.}}%
}
\makeatother

\begin{document}

\subjclass{Primary 11R23, 37A05, 37A44. Secondary 11K06, 11K41}
\keywords{Iwasawa theory, Ergodic theory, Bernoulli shifts, equidistribution modulo $1$}

\maketitle
\begin{abstract} 
An important ingredient in the Ferrero--Washington proof of the vanishing of cyclotomic $\mu$-invariant for Kubota--Leopoldt $p$-adic $L$-functions is an equidistribution result which they established using the Weyl criterion. The purpose of our manuscript is to provide an alternative proof by adopting a dynamical approach. A key ingredient to our methods is studying an ergodic skew-product map on $\Z_p \times [0,1]$, which is then suitably identified as a factor  of the $2$-sided Bernoulli shift on the sample space $\{0,1,2,\cdots,p-1\}^{\mathbb{Z}}$. \end{abstract}

\tableofcontents
 
\section{Introduction}
Let $p$ denote an odd prime. Let $p^{e_n}$ denote the exact power of $p$ dividing the class number of $\Q(\mu_{p^{n+1}})$. Iwasawa \cite{MR124320} has showed that there exists constants $\lambda$, $\mu$ and $\nu$ (depending only on $p$) such that  
\[e_n = \lambda n + \mu p^n + \nu, \text{ for all sufficiently large } n.\]
Iwasawa further conjectured that the invariant $\mu$ equals zero. A celebrated result of Ferrero--Washington \cite{MR528968} established this conjecture. 

\begin{Theorem}[Ferrero--Washington]\label{thm:FW}
$\mu=0$.	
\end{Theorem}

Before detailing the motivations behind our present work and to better place it in context, it will be helpful to first briefly highlight some of the main ingredients involved in the proof by Ferrero--Washington.  Let us introduce a few notations for this purpose. For every $\alpha$ in $\Z_p$, we can consider its $p$-adic expansion:
\begin{align}\label{eq:padicexpansion}
\alpha &= t_0(\alpha) + t_1(\alpha)p^1 + t_2(\alpha)p^2 + \cdots +t_n(\alpha)p^n + t_{n+1}(\alpha)p^{n+1} +  \cdots,
\end{align}
along with its associated partial sums for each $n \geq 1$:
\begin{align}
	s_{n}(\alpha) = t_0(\alpha) + t_1(\alpha)p^1 + t_2(\alpha)p^2 + \cdots + t_{n}(\alpha)p^{n}.
\end{align}
Here, the digits $t_n(\alpha)$'s belong to the set $\{0,1,2,\cdots, p-1\}$. Notice that the element $p^{-n-1}s_{n}(\alpha)$ always belongs to the interval $[0,1)$. \\

The following proposition is the main criterion (originally due to Iwasawa \cite{MR0124317}) utilized by Ferrero--Washington in their proof of Theorem \ref{thm:FW}.  

\begin{Proposition}[Iwasawa, Ferrero--Washington]\label{prop:criterion} The following statements are equivalent:
\begin{enumerate}
\item $\mu >0$.
\item There exists an odd integer $3 \leq d \leq p-2$ such that for all $n \geq 0$ and $\alpha \in \Z_p$, we have
\begin{align}\label{eq:criterion}
\sum_{\eta^{p-1}=1} t_{n}(\alpha\eta)\eta^d \equiv 0 \ (\mathrm{mod} \ p).	
\end{align}
\end{enumerate}
\end{Proposition}

There are two main ingredients used to establish this proposition. The first main ingredient is an application of the reflection theorem (Spiegelungssatz); this ingredient is used to obtain the restriction to odd integers $3 \leq d \leq p-2$. The second main ingredient is an application of the Stickelberger theorem; one uses the fact that the Stickelberger element 

\begin{align}
\sum_{\substack{u=1 \\ u \equiv 1 \ \mathrm{mod} \ p}}^{p^{n+1}} \left(\sum_{\eta^{p-1}=1} \dfrac{s_{n}(u \eta)}{p^{n+1}} \ \sigma_\eta^{-1} \right) \sigma_{u}^{-1} \qquad \in \  \Q[(\Z/p^{n+1})^\times]
\end{align}

annihilates the minus part of the class group of $\Q(\mu_{p^{n+1}})$.  Here, $\sigma_u$ and $\sigma_\eta$ are the group ring elements in the group ring $\Q[(\Z/p^{n+1})^\times]$ corresponding to $u$ and $\eta$ in the group $(\Z/p^{n+1})^\times$. The connection with $p$-adic $L$-functions arises since a compatible family of Stickelberger elements can be used to construct the Kubota--Leopoldt $p$-adic $L$-function (see Iwasawa's book \cite{MR0360526}). The explicit form of the Stickelberger element is what explains the presence of the digits $t_n$'s in the criterion given in equation (\ref{eq:criterion}). \\

To establish $\mu=0$, Ferrero--Washington have to show that equation (\ref{eq:criterion}) cannot hold. To do so,  Ferrero--Washington use the following proposition. It is the key equidistribution result required by them. It is also the equidistribution result alluded to in the title of our manuscript.

\begin{Proposition}[Ferrero--Washington]\label{prop:equidistribution}
Let $r \geq 1$ be an integer. Let $\{\beta_1,\cdots, \beta_r\} \subset \Z_p$ be a linearly independent set over $\Q$. Consider the following subset of $\Z_p$: 
\begin{align}\label{eq:equidistributedseq}
G_r=\left\{
	\begin{array}{ll} 
    	\alpha \in \Z_p,  & \text{ such that the sequence }  \left\{\bigg(p^{-n-1}s_{n}(\alpha\beta_1),\cdots,p^{-n-1}s_{n}(\alpha\beta_r)\bigg)\right\}_{n=0}^{\infty} 
    	 \text{ is } \\ & \text{ equidistributed in }  [0,1]^r \text{ with respect to the standard Borel measure.}
    	 \end{array} \right\}.	
\end{align}
Then, $G_r$ has full Haar measure in $\Z_p$.
\end{Proposition}

The proof of Proposition \ref{prop:equidistribution} by Ferrero--Washington uses the Weyl criterion. Our objective is to reprove this proposition from a dynamical perspective by utilizing the following fact: almost every point on a compact space $X$ is generic with respect to a continuous ergodic map $T$ on $X$ and a measure on $X$ that is $T$-invariant. See \cite[Corollary 4.20]{MR2723325}. The space $X$ in our set up will turn out to be an $r$-fold product of $\{0,1,2,\cdots,p-1\}^{\mathbb{Z}}$ equipped with the uniform probability measure, while the ergodic map $T$ will turn out to be an $r$-fold product of the two-sided Bernoulli shift. What we mean by a dynamical perspective is that the equidistributed sequence in equation (\ref{eq:equidistributedseq}) will be obtained from the orbit $\{T^{[n]}(x)\}$ of a generic point for the action of the $2$-sided Bernoulli shift $T$ on $X$. We expand on this point in Sections \ref{sec:intro_bernoulli} and \ref{sec:proof}.
\begin{remark}
The equidistribution statement in equation (\ref{eq:equidistributedseq}) can also equivalently be stated in terms of the standard Lebesgue measure as it is the completion of the standard Borel measure.  	
\end{remark}

\begin{remark}
Our main focus in this manuscript is the equidistribution result of Ferrero--Washington \cite{MR528968}. We explain our motivations in the next subsection. We simply refer the interested reader to their original article on how one can deduce $\mu=0$ from the equidistribution result. Ferrero--Washington crucially use \cite[Proposition 3]{MR528968} for this deduction step. Although the implication from the equidistribution result to the $\mu=0$ result is not ergodic in essence, it nevertheless appears to be non-trivial.  In fact, there is a formal resemblance between \cite[Proposition 3]{MR528968} and \cite[Proposition 3.1]{MR0732547} that was used by Sinnott in his independent proof of $\mu=0$.  See \cite[Appendix]{MR1097627}.
\end{remark}

\subsection{Motivations}

There are historical precedents to providing alternate proofs of equidistribution results using approaches from ergodic theory. One very well-known example of such a result is Furstenberg's ergodic-theoretic proof (see \cite{MR603625}) of the equidistribution in $[0,1]$ of the sequence $\{p(n)\}_{n=1}^\infty$, for any polynomial $p(x)$ with an irrational leading coefficient. It was the ergodic theoretic approach that led to Furstenberg's famous $\times p,\times q$ conjectures. Another inspiration arises from works of Cornut and Vatsal \cite{MR1908058,MR1892842,MR1953292} who successfully applied Ratner's theorems to prove results in anticyclotomic Iwasawa theory. Our first motivation is to follow these precedents. \\

One may reasonably question the necessity to provide an alternate proof of Proposition \ref{prop:equidistribution} since the original proof is a simple one using the Weyl criterion. The application of the Weyl criterion in \cite{MR528968}, however, relies on computing some explicit (although fairly elementary) bounds using products of exponential test functions. Our second motivation in considering a dynamic set up is to avoid this explicit analysis  so that the method is more amenable to generalizations (where the test functions might be more complicated to compute explicit bounds). One concrete generalization to consider would be quotients  of $\SL_2(\mathbb{R})$ by Fuchsian groups (such as finite index subgroups of $\SL_2(\Z)$).   \\

As we explain in Section \ref{sec:intro_bernoulli}, the connection with Bernoulli shifts stems from two observations: (i) the space $\{0,1,2,\cdots,p-1\}^{\mathbb{Z}}$ is measurably isomorphic to 
\begin{align}\label{eq:solenoid}
\Z[1/p] \backslash \left({\Q_p \times \R}\right).
\end{align}
(ii) Under this identification, the shift map corresponds to multiplication by $1/p$ map. The space in equation (\ref{eq:solenoid}) is the so-called ``\textit{$p$-adic solenoid}''. The $\times 1/p$ map could be viewed as an analog of the ``\textit{Hecke map}'', e.g. consider  on $\SL_2(\Z[1/p]) \backslash \left({\SL_2(\Q_p) \times \SL_2(\R)}\right)$, the $p$-adic extension of $\SL_2(\Z)/\SL_2(\R)$, the right-multiplication by $\left[\begin{array}{cc}p^{-1} & 0 \\ 0 & p\end{array}\right]$.

The advantage of using a dynamical model is that we have more ergodic tools at our disposal.  It will be interesting to translate these ergodic properties to new concrete results in Iwasawa theory.  For example, in the set up involving quotients of $\SL_2(\mathbb{R})$ by Fuchsian groups, we have the ergodicity of geodesic and unipotent flows at our disposal. In this case, there is a concrete connection to symbolic dynamics via the Bowen--Series coding \cite{MR556585}. Beginning with the connection between geodesic flows and continued fractions, this topic has a long and rich history following numerous authors such as Adler--Flatto, Artin, Bowen, Hedlund, Hopf,  Katok, Martin, Manning, Series, etc. Our third motivation, in our setting, lies in identifying the connection between Bernoulli shifts (more generally, symbolic dynamics) and Iwasawa theory. The theory of Bernoulli shifts is classical and has been very well-studied. For this reason, we hope that this connection will be both fruitful and of independent interest.   \\

Our approach involving ergodic theory is rooted in the principle of transferring dynamics from the real to the $p$-adic side. This principle has been applied with great success in ergodic theory. For an application of this principle in a setting similar to ours, see, for example, Lindenstrauss's work \cite{MR1826489}. In our situation, this principle allows us to concretely reduce from the general $r \geq 1$ case to the $r=1$ case. See Propositions \ref{prop:nodependenceonrealcoordinates} and  \ref{prop:reductionstep}. This reduction step suggests (at least to us) that the $r=1$ case entirely encapsulates the ergodic nature of the problem at hand. Our fourth motivation involves emphasizing the $r=1$ case, which seems important, especially in light of conjectures of Mazur--Rubin \cite{mazur2021arithmetic} on residual equidistribution of rational modular symbols. See \cite[Conjecture B(3)]{lee2019dynamics}.

One can compare this with \cite[Section 2.1]{MR2275607}, where a similar reduction step, in the spirit of ergodic theory, is suggested via analogy with Kronecker's classical result by considering maps from the $1$-dimensional space to the general $r$-dimensional space. However, for the precise statement involving full Haar measure in Proposition \ref{prop:equidistribution}, it seems to us that we would need to instead consider maps from the general $r$-dimensional space to the $1$-dimensional space. This reduction step, seems to us, to be group theoretic (or representation theoretic) in its essence, rather than involving genuine inputs from ergodic theory.

\subsection{The $p$-adic solenoid and Bernoulli shifts} \label{sec:intro_bernoulli}

In our setting, there are five spaces (labeled $X_0$, $X_1$, $X_2$, $X_3$ and $X_4$) on which we can consider dynamics. These five spaces are all measurably isomorphic; so choosing which space to work with boils down to the perspective the reader is most comfortable with. 

\begin{itemize} 
\item The first space $X_1$ arises as a $p$-adic extension to the cirlce $\Z \backslash \R$. We have the following natural inclusion map:
 \begin{align} \label{eq:inclusionlattice}
 \Z[1/p] &\hookrightarrow \left(\Q_p \times \R\right), \\
 x &\mapsto  (-x,x). \notag
 \end{align}
Endow $\R$ and $\Q_p$ with the usual archimedean and $p$-adic topologies respectively. Under the inclusion given above, $\Z[1/p]$ sits inside $\Q_p \times \R$ as a discrete lattice. We will be interested in the quotient space 
\[X_1 \coloneqq \Z[1/p] \backslash \left(\Q_p \times \R\right),\]
which is also often called the $p$-adic solenoid.  See \cite[Appendix to Chapter 1]{MR1760253} for more details on the $p$-adic solenoid. The space $X_1$ is compact. We let $\nu_1$ be the Haar measure on $X_1$, normalized so that it is a probability measure. We have the $\times
1/p$-map on $X_1$ given below:
\begin{align}
T_1 : X_1 & \rightarrow X_1, \\
(\alpha,x) & \mapsto \left(\dfrac{\alpha}{p},\dfrac{x}{p}\right). \notag
\end{align}

The self-map $T_1$ is measurable with respect to $\nu_1$.  

\item Under the inclusion given in equation (\ref{eq:inclusionlattice}), we see that $\Z$ sits as a discrete lattice inside $\Z_p \times \R$. One can also consider the corresponding quotient space 
\[X_2 \coloneqq \Z \backslash \left(\Z_p \times \R\right).\]
For every $\beta$ in $\Q_p$, there exists an unique element in $\Z_p$, denoted $\lfloor \beta \rfloor$,  such that $\beta - \lfloor \beta \rfloor$ belongs to $\Q$ and satisfies
\[0 \leq \beta - \lfloor \beta \rfloor  <1.\] 
Using this observation, it is not hard to see that the natural inclusion $\Z_p \times \R \hookrightarrow \Q_p \times \R$ induces the following isomorphism of topological spaces: 
\[X_2 \cong X_1.\]
Let $\nu_2$ be the induced measure on $X_2$ under this isomorphism, which coincides with the Haar measure normalised to be a probability measure. Under this isomorphism, the self-map on $X_2$ induced by $T_1$ can be described as follows:
\begin{align}\label{eq:formulashiftZp}
T_2 : X_2 &\rightarrow X_2, \\
	 \left(\alpha,x\right) & \mapsto \left(\dfrac{\alpha - t_0(\alpha)}{p},\dfrac{x+t_0(\alpha)}{p}\right) \notag
\end{align}
The self-map $T_2$ is measurable with respect to $\nu_2$.

\item We will consider the inverse limit $X_0 \coloneqq \varprojlim_n \dfrac{\R}{p^n\Z}$ of topological groups, given by the natural surjections $\dfrac{\R}{p^n\Z} \twoheadrightarrow \dfrac{\R}{p^{n-1}\Z}$. We have the following natural isomorphism of topological groups (see \cite[Appendix to Chapter 1]{MR1760253}): 
\begin{align*}
	X_2 & \xrightarrow {\cong} X_0, \\
	(\beta, x) & \mapsto \left(x, \  x+ t_0(\beta), \ \cdots, \  x+t_0(\beta) + \cdots t_{n-1}(\beta)p^{n-1}, \cdots, \ \right)
\end{align*}

The measure $\nu_0$ on $X_0$ will be the pushforward of the measure $\nu_2$ under this isomorphism.  The map $T_0: X_0 \rightarrow X_0$, associated via this isomorphism, to $T_2$ is given below:
\begin{align*}
	\left(x, \  x+ t_0(\beta), \ \cdots, \  x+t_0(\beta) + \cdots t_{n-1}(\beta)p^{n-1}, \cdots, \ \right)  \longmapsto \\     \left(\dfrac{x+ t_0(\beta)}{p}, \ \cdots, \  \dfrac{x+t_0(\beta) + \cdots t_{n-1}(\beta)p^{n-1}}{p}, \cdots, \ \right)  \end{align*}

\item We also consider the space 
\begin{align}
\qquad X_3 \coloneqq  \Z_p \times [0,1].	
\end{align}
as a choice of a ``\textit{fundamental domain}'' for $X_1$ and $X_2$. We equip $X_3$ with the measure $\nu_{\mathrm{Bor}} \times \nu_{\mathrm{Haar}}$, which we denote $\nu_3$. Here,$\nu_{\mathrm{Bor}}$ is the standard Borel measure on $[0,1]$, while $\nu_{\mathrm{Haar}}$ is the usual Haar measure on $\Z_p$. While the descriptions of $X_1$ and $X_2$ arise naturally in theory, $X_3$ will turn out to be the space that is the most transparent to work with.

There is a natural continuous surjection
\begin{align}\label{eq:surjectionX3X2}
	X_3 \twoheadrightarrow X_2.
\end{align}
It is straightforward to check that under this surjection, the pushforward measure of $\nu_3$ on $X_2$ coincides with the measure $\nu_2$. Furthermore, it is also straightforward to check that this surjection is an isomorphism outside a set of measure zero. We have a measurable self-map \[T_3: X_3\rightarrow X_3,\] which is given by the same formula as in equation (\ref{eq:formulashiftZp}) and hence is compatible with $T_2$ and the surjection in equation (\ref{eq:surjectionX3X2}). This map is often referred to as a \textit{skew-product} map.

\item Lastly, we have the following space on which the two-sided Bernoulli shift acts:
\begin{align}
 X_4 \coloneqq \{0,1,2,\cdots,p-1\}^{\mathbb{Z}}.
 \end{align} 
Every element in $X_4$ is of the form $\left(\cdots, b_{-n},\cdots, b_{-1} \vert b_0,b_1,\cdots,b_n,\cdots\right)$. We view $X_4$ as a Bernoulli scheme (as in \cite[Example 2.9]{MR2723325}), equipped with the probability measure (denoted $\nu_4)$ corresponding to the uniform probability vector $(1/p,1/p,\cdots,1/p)$ on the sample space $\{0,1,\cdots,p-1\}$. 

For every $x$ in $[0,1]$, we can consider an expansion in base $p$:
\begin{align}
	x = \dfrac{a_{-1}(x)}{p} + \dfrac{a_{-2}(x)}{p^2} + \cdots \dfrac{a_{-n}(x)}{p^n} + \cdots
\end{align}
Here, each of the $a_i(x)$'s belong to the set $\{0,1,\cdots,p-1\}$. We have a natural continuous surjection 
\begin{align*}
X_4 & \rightarrow X_3,\\
\left( \cdots, a_{-n}, \cdots, a_{-2}, a_{-1} \ \biggl\vert \  t_0, t_1, t_2, \cdots, t_n,\cdots\right) &\mapsto \left(\sum_{i=0}^\infty t_i p^i ,\ \sum \limits_{j=0}^
\infty \dfrac{a_{-(j+1)}}{p^{j+1}}\right).
\end{align*}

This map $X_4 \rightarrow X_3$ is essentially obtained from the map $\Z_p \rightarrow [0,1]$ sending a $p$-adic expansion to a $p$-ary expansion. It is this map that furnishes the connection to symbolic dynamics. For this reason, although it is fairly elementary and may be well-known to experts, we show in Section \ref{sec:compatibility} that under the $X_4\rightarrow X_3$, the pushforward of $\nu_4$ on $X_3$ coincides with the measure $\nu_3$. This map turns out to be a bijection outside a set of measure zero (corresponding to the rational numbers whose denominators are powers of $p$). We can consider the following measurable self-map, which is often called the Bernoulli shift map:
\begin{align}
T_4: X_4 &\rightarrow X_4,\\ 
\left(\cdots, b_{-n},\cdots, b_{-1} \vert b_0,b_1,\cdots,b_n,\cdots\right) &\mapsto \left(\cdots, b_{-n},\cdots,  b_{-1}, b_0 \vert  b_1,\cdots,b_n,\cdots\right). \notag
\end{align}
The maps $T_3$ and $T_4$ are compatible with the surjection $X_4 \twoheadrightarrow X_3$. The map $T_4$ is ergodic with respect to $\nu_4$. See \cite[Proposition 2.15]{MR2723325}. In fact, the Bernoulli-shift map $T_4$ is mixing and as a result, for each $r \geq 1$, the $r$-fold product $T_4^r$ is also ergodic. See \cite[Exercise 2.7.9 and Theorem 2.36]{MR2723325}.

 \end{itemize}

\begin{figure}
\begin{flushright}
	\begin{tikzpicture}
		\draw[thick](0,0) circle (1.5);
				\draw[thick](5,0) circle (1.5);
				\draw[thick](-5,0) circle (1.5);
			  \draw[very thick,->] (2,0) -- (3,0); 	 
			  \draw[very thick,->] (-3,0) -- (-2,0); 
			    \draw[very thick,->] (-8,0) -- (-7,0); 
			       \node(a1) at (2.5,0.5) {$p$}; 
			         \node(a2) at (-2.5,0.5) {$p$}; 
			           \node(b) at (-7.5,0.5) {$p$}; 
			       \node[circle,fill=black, inner sep=0pt,minimum size=2pt] (b1) at (-9,0) {};  
 \node[circle,fill=black, inner sep=0pt,minimum size=2pt] (b2) at (-8.4,0) {};  
 \node[circle,fill=black, inner sep=0pt,minimum size=2pt] (b3) at (-8.7,0) {}; 
 

   \draw [thick,domain=-180:180,fill=spanishgray] plot ({6.5+\radiusc*cos(\x)}, {\radiusc*sin(\x)});

   \draw [thick,domain=-180:180,fill=spanishgray] plot ({1.5+\radiusc*cos(\x)}, {\radiusc*sin(\x)});

  \draw [thick,domain=-180:180,fill=darkslategray] plot ({1.5*cos(120)+\radiusc*cos(\x)}, {1.5*sin(120)+\radiusc*sin(\x)});

  \draw [thick,domain=-180:180,fill=white] plot ({1.5*cos(240)+\radiusc*cos(\x)}, {1.5*sin(240)+\radiusc*sin(\x)});


 \draw [thick,domain=-180:180,fill=spanishgray]plot ({(-5+1.5*cos(\anglezerozero))+\radiusc*cos(\x)}, {(0+1.5*sin(\anglezerozero))+\radiusc*sin(\x)}); 
\draw [thick,domain=-180:180,fill=spanishgray]plot ({(-5+1.5*cos(\anglezerozero))+\radiusb*cos(\x)}, {(0+1.5*sin(\anglezerozero))+\radiusb*sin(\x)});

\draw [thick,domain=-180:180,fill=spanishgray]plot ({(-5+1.5*cos(\anglezeroone))+\radiusc*cos(\x)}, {(0+1.5*sin(\anglezeroone))+\radiusc*sin(\x)}); 
\draw [thick,domain=-180:180,fill=darkslategray]plot ({(-5+1.5*cos(\anglezeroone))+\radiusb*cos(\x)}, {(0+1.5*sin(\anglezeroone))+\radiusb*sin(\x)});

\draw [thick,domain=-180:180,fill=spanishgray]plot ({(-5+1.5*cos(\anglezerotwo))+\radiusc*cos(\x)}, {(0+1.5*sin(\anglezerotwo))+\radiusc*sin(\x)}); 
\draw [thick,domain=-180:180,fill=white]plot ({(-5+1.5*cos(\anglezerotwo))+\radiusb*cos(\x)}, {(0+1.5*sin(\anglezerotwo))+\radiusb*sin(\x)});

  \draw [thick,domain=-180:180,fill=darkslategray]plot ({(-5+1.5*cos(\angleonezero))+\radiusc*cos(\x)}, {(0+1.5*sin(\angleonezero))+\radiusc*sin(\x)}); 
\draw [thick,domain=-180:180,fill=spanishgray]plot ({(-5+1.5*cos(\angleonezero))+\radiusb*cos(\x)}, {(0+1.5*sin(\angleonezero))+\radiusb*sin(\x)});

   \draw [thick,domain=-180:180,fill=darkslategray]plot ({(-5+1.5*cos(\angleoneone))+\radiusc*cos(\x)}, {(0+1.5*sin(\angleoneone))+\radiusc*sin(\x)}); 
\draw [thick,domain=-180:180,fill=darkslategray]plot ({(-5+1.5*cos(\angleoneone))+\radiusb*cos(\x)}, {(0+1.5*sin(\angleoneone))+\radiusb*sin(\x)});

   \draw [thick,domain=-180:180,fill=darkslategray]plot ({(-5+1.5*cos(\angleonetwo))+\radiusc*cos(\x)}, {(0+1.5*sin(\angleonetwo))+\radiusc*sin(\x)}); 
\draw [thick,domain=-180:180,fill=white]plot ({(-5+1.5*cos(\angleonetwo))+\radiusb*cos(\x)}, {(0+1.5*sin(\angleonetwo))+\radiusb*sin(\x)});

   \draw [thick,domain=-180:180,fill=white]plot ({(-5+1.5*cos(\angletwozero))+\radiusc*cos(\x)}, {(0+1.5*sin(\angletwozero))+\radiusc*sin(\x)}); 
\draw [thick,domain=-180:180,fill=spanishgray]plot ({(-5+1.5*cos(\angletwozero))+\radiusb*cos(\x)}, {(0+1.5*sin(\angletwozero))+\radiusb*sin(\x)});

   \draw [thick,domain=-180:180,fill=white]plot ({(-5+1.5*cos(\angletwoone))+\radiusc*cos(\x)}, {(0+1.5*sin(\angletwoone))+\radiusc*sin(\x)}); 
\draw [thick,domain=-180:180,fill=darkslategray]plot ({(-5+1.5*cos(\angletwoone))+\radiusb*cos(\x)}, {(0+1.5*sin(\angletwoone))+\radiusb*sin(\x)}); 

   \draw [thick,domain=-180:180,fill=white]plot ({(-5+1.5*cos(\angletwotwo))+\radiusc*cos(\x)}, {(0+1.5*sin(\angletwotwo))+\radiusc*sin(\x)}); 
\draw [thick,domain=-180:180,fill=white]plot ({(-5+1.5*cos(\angletwotwo))+\radiusb*cos(\x)}, {(0+1.5*sin(\angletwotwo))+\radiusb*sin(\x)});   
 
	\end{tikzpicture}
	\caption[]
	{
	\begin{minipage}{0.85\textwidth}
{The $p$-adic solenoid can be visualized as the inverse limit $\cdots \xrightarrow {p} \frac{\mathbb{R}}{\mathbb{Z}} \xrightarrow{p} \frac{\mathbb{R}}{\mathbb{Z}}\xrightarrow{p}  \frac{\mathbb{R}}{\mathbb{Z}}$. \\ In this figure, we consider $p=3$. The color {\color{spanishgray} gray}, {\color{darkslategray} dark slate gray} and white represent the digits $0$, $1$ and $2$ respectively. Here, the unit circle is used to represent $\R/\Z$. Here, the node \begin{tikzpicture}
	   \draw [thick,domain=-180:180,fill=spanishgray]plot ({(0))+\radiusc*cos(\x)}, {(0)+\radiusc*sin(\x)}); 
\draw [thick,domain=-180:180,fill=darkslategray]plot ({(0))+\radiusb*cos(\x)}, {(0)+\radiusb*sin(\x)}); 
\end{tikzpicture} represents the point $\exp\left(2\pi\times \dfrac{({\color{spanishgray}0}+{\color{darkslategray}1}\times 3)}{3^2}\right)$ on the unit circle.}
\end{minipage}}
\end{flushright}
\end{figure}

To summarize, we have the following commutative diagram of measurable spaces:
 \begin{align*}
 \xymatrix{
 (X_4,\nu_4) \ar@{->>}[r]\ar[d]^{T_4}&  (X_3,\nu_3) \ar@{->>}[r]\ar[d]^{T_3}&  (X_2,\nu_2) \ar[r]^{\cong}\ar[d]^{T_2}&  (X_1,\nu_1) \ar[d]^{T_1}  \ar@{->>}[r]^{\cong}& (X_0,\nu_0) \ar[d]^{T_0} \\
  (X_4,\nu_4) \ar@{->>}[r] &  (X_3,\nu_3) \ar@{->>}[r]&  (X_2,\nu_2)\ar[r]^{\cong}& (X_1,\nu_1) \ar@{->>}[r]^{\cong}& (X_0,\nu_0)\\
 }	
 \end{align*}

Let $r \geq 1$ be an integer. For each $i \in \left\{0,1,2,3,4\right\}$, we can consider the product measure $\nu_i^r$ on $X_i^r$ along with the natural extensions to the maps $T_i$ as follows:
\begin{align*}
T_i^r : X_i^r &\rightarrow X_i^r, \\ 
	\left(x_1,\cdots x_r\right) &\mapsto \left(T_i(x_1),\cdots, T_i(x_r)\right).
\end{align*}

The self-map $T_i^r$ is measurable with respect to $\nu_i^r$.  We have the following commutative diagram of measurable spaces:
\begin{align*}
 \xymatrix{
 (X_4^r,\nu_4^r) \ar@{->>}[r] \ar[d]& (X_3^r,\nu_3^r)  \ar@{->>}[r]\ar[d]& (X_2^r,\nu_2^r) \ar[r]^{\cong}\ar[d]& (X_1^r,\nu_1^r) \ar[r]^{\cong}\ar[d] &  (X_0^r,\nu_0^r) \ar[d]\\
  (X_4^r,\nu_4^r) \ar@{->>}[r] &(X_3^r,\nu_3^r)  \ar@{->>}[r]& (X_2^r,\nu_2^r)  \ar[r]^{\cong}& (X_1^r,\nu_1^r) \ar[r]^{\cong}& (X_0^r,\nu_0^r)  \\
 }	
 \end{align*}

Since the spaces $X_i^r$ are all measurably isomorphic to each other and since $T_4^r$ is ergodic, each of the maps $T_i^r$ is ergodic. See, for example, \cite[Exercise 2.3.4]{MR2723325}. We also have the following related self-map which is measurable with respect to the Haar measure on $\Z_p^r$:
\begin{align}
	\mathcal{T}_3^r: \Z_p^r &\rightarrow \Z_p^r, \\
	\left(\alpha_1,\cdots,\alpha_r\right) &\mapsto \left(\dfrac{\alpha_1 - t_0(\alpha_1)}{p},\cdots,\dfrac{\alpha_r - t_0(\alpha_r)}{p} \right). \notag
\end{align}

\begin{remark}
Given a triple $(X,T,\nu)$, a point $x \in X$ is said to be generic for the map $T$ if its orbit, given by the set $\{x, T(x), T^2(x), \ldots,\}$, is equidistributed with respect to the measure $\nu$. In \cite{MR528968}, a generic point $(\gamma_1,\cdots,\gamma_r)$ for the map $\mathcal{T}_3^r$ is called \textit{jointly normal}. A generic point $\gamma$ for the map $\mathcal{T}_3$ is also simply called \textit{normal}.	 
\end{remark}
\begin{remark}\label{rem:1sidedbernoullishift}
Just as the map $T_3$ is analogous to the $2$-sided Bernoulli shift, the map $\mathcal{T}_3$ is analogous to the $1$-sided Bernoulli shift on $\{0,1,\cdots,p-1\}$. That is, $(\Z_p,\mathcal{T}_3,\nu_{\mathrm{Haar}})$ is measurably isomorphic to the $1$-sided Bernoulli shift on $\left\{0,1,\cdots,p-1\right\}^{\mathbb{N}}$ equipped with the probability measure induced by the uniform probability vector $\left(1/p,1/p,\cdots, 1/p\right)$ on the sample space $\{0,1,\cdots,p-1\}$. In fact, we have a homeomorphism of topological spaces. The one-sided Bernoulli shift map and its $r$-fold products are also ergodic. See Exercise 2.3.4 along with  Proposition 2.15 in \cite{MR2723325}. 
\end{remark}

\subsection{Proof of the equidistribution result}\label{sec:proof}
In this subsection, we show how Proposition \ref{prop:equidistribution} follows from Propositions \ref{prop:nodependenceonrealcoordinates} and \ref{prop:reductionstep}. The proofs of Propositions \ref{prop:nodependenceonrealcoordinates} and \ref{prop:reductionstep} are later given in Sections \ref{sec:proofpropnodependence} and \ref{sec:proofpropreduction}. We follow the notations of Proposition \ref{prop:equidistribution}.

\begin{proposition}\label{prop:nodependenceonrealcoordinates}
The following statements are equivalent:
\begin{enumerate}
\item\label{it:prop1} For every $(x_1,\cdots,x_r)$ in $[0,1]^r$, the element $\big(\left(\gamma_1,x_1\right),\left(\gamma_2,x_2\right),\cdots,\left(\gamma_r,x_r\right)\big)$ in $X_3^r$ is a generic point for $T_3^r$.
\item\label{it:prop2}   $\big(\left(\gamma_1,0\right),\left(\gamma_2,0\right),\cdots,\left(\gamma_r,0\right)\big)$ in $X_3^r$ is a generic point for $T_3^r$.
\item\label{it:prop3}  $(\gamma_1,\cdots,\gamma_r)$ is a generic point for $\mathcal{T}_3^r$.
\item \label{it:prop4} For every $(x_1,\cdots,x_r)$ in $\mathbb{R}^r$, the element $\big(\left(\gamma_1,x_1\right),\left(\gamma_2,x_2\right),\cdots,\left(\gamma_r,x_r\right)\big)$ in $X_2^r$ is a generic point for $T_2^r$. 
\item \label{it:prop5} For every $(x_1,\cdots,x_r)$ in  $\mathbb{R}^r$, the element $\big(\left(\gamma_1,x_1\right),\left(\gamma_2,x_2\right),\cdots,\left(\gamma_r,x_r\right)\big)$ in $X_1^r$ is a generic point for $T_1^r$.
\item \label{it:prop6} $\big(\left(\gamma_1,0\right),\left(\gamma_2,0\right),\cdots,\left(\gamma_r,0\right)\big)$ in $X_2^r$ is a generic point for $T_2^r$.
\item \label{it:prop7} $\big(\left(\gamma_1,0\right),\left(\gamma_2,0\right),\cdots,\left(\gamma_r,0\right)\big)$ in $X_1^r$ is a generic point for $T_1^r$.
\end{enumerate}

\end{proposition}
Throughout, we assume that $\{\beta_1,\cdots,\beta_r\}$ is a linearly independent set. Consider the following countable subset in $\Z_p$: 
\begin{align}\label{eq:discretesubset}
V \coloneqq \left\{m_1 \beta_1 + \cdots m_r \beta_r \in \Z_p, \text{ such that } (m_1,\cdots,m_r) \in \Z^r \setminus (0,\cdots,0) \right\}.	
\end{align}

Since $\{\beta_1,\cdots,\beta_r\}$ is a linearly independent set, observe that for each $(m_1,\cdots,m_r) \in \Z^r \setminus (0,\cdots,0)$, the linear combination $m_1 \beta_1 + \cdots m_r \beta_r$ is never equal to $0$. For each such a linear combination $\sigma$, we can write $\sigma$ as $p^{\mathrm{val}(\sigma)}u_\sigma$, for a unique element $u_\sigma$ in $\Z_p^\times$ and a non-negative integer $\mathrm{val}(\sigma)$.

\begin{proposition}\label{prop:reductionstep}
Let $\alpha$ be an element of $\Z_p$. The following statements are equivalent.
\begin{enumerate}[label=(\arabic*),ref=(\arabic*)]
\item\label{it:propreduction1} $\left(\alpha\beta_1,\cdots,\alpha\beta_r\right)$ is a generic point for $\mathcal{T}_3^r$.
\item\label{it:propreduction2} For every $\sigma$ in $V$, the element $\sigma\alpha$ is a generic point for $\mathcal{T}_3$.	
\end{enumerate}
\end{proposition}

\subsubsection{Propositions \ref{prop:nodependenceonrealcoordinates} and \ref{prop:reductionstep} $\implies$ Proposition \ref{prop:equidistribution}} \label{subsec:proofimplication}

We now proceed to prove Proposition \ref{prop:equidistribution} assuming that Propositions \ref{prop:nodependenceonrealcoordinates} and \ref{prop:reductionstep} hold. Consider the following set:
\begin{align}
H_r =\left\{\alpha \in \Z_p, \text{ such that } \big(\left(\alpha\beta_1,0\right),\cdots,\left(\alpha\beta_r,0\right)\big) \text{ is a generic point for } T_3^r\right\}.
\end{align}

If $n\geq 0$ is an integer, observe first that 
\begin{align*}
	& (T_3^r)^{[n+1]}\big(\left(\gamma_1,0\right),\left(\gamma_2,0\right),\cdots,\left(\gamma_r,0\right)\big) = \\  & \qquad \qquad \left(\left(\dfrac{\gamma_1-s_n(\gamma_1)}{p^{n+1}},\dfrac{s_n(\gamma_1)}{p^{n+1}}\right),\left(\dfrac{\gamma_2-s_n(\gamma_2)}{p^{n+1}},\dfrac{s_n(\gamma_2)}{p^{n+1}}\right), \cdots, \left(\dfrac{\gamma_r-s_n(\gamma_r)}{p^{n+1}},\dfrac{s_n(\gamma_r)}{p^{n+1}}\right)\right).
\end{align*}
If $\big(\left(\gamma_1,0\right),\left(\gamma_2,0\right),\cdots,\left(\gamma_r,0\right)\big)$ is a generic point in $\big(\Z_p \times [0,1])^r$ for the map $T_3^r$, then the sequence \[\bigg\{(T_3^r)^{[n+1]}\big(\left(\gamma_1,0\right),\left(\gamma_2,0\right),\cdots,\left(\gamma_r,0\right)\big)\bigg\}_{n=1}^{\infty}\] is equidistributed in $\big(\Z_p \times [0,1]\big)^r$ with respect to $\nu_3^r$. Let $\varpi$ denote the natural continuous surjection $\big(\Z_p \times [0,1]\big)^r \rightarrow [0,1]^r$. Since the measure $\nu_3^r$ equals $(\nu_{\mathrm{Haar}} \times \nu_{\mathrm{Bor}})^r$ on $(\Z_p\times[0,1])^r$, the pushforward measure $\varpi^*(\nu_3^r)$ coincides with $\nu_{\mathrm{Bor}}^r$. As a result, for any continuous function $f: [0,1]^r \rightarrow \R$, we have (see, for instance \cite[Theorem 3.6.1]{MR2267655}, a standard result relating integrals with respect to pushforward measures)
\[ \int (f  \circ \pi)\ d\nu_3^r=  \int f \ d\nu_{\mathrm{Bor}}^r.\]

These observations let us conclude that the sequence
\[\bigg\{\big(p^{-n-1}s_n(\gamma_1),\ p^{-n-1}s_n(\gamma_2),\cdots,p^{-n-1}s_n(\gamma_r)\big)\bigg\}_{n=1}^\infty\]
must then be equidistributed in $[0,1]^r$ with respect to the Borel measure. Using this observation, one can conclude that 
\begin{align}
H_r \subseteq G_r.	
\end{align}

To prove Proposition \ref{prop:equidistribution}, it is enough to show that $H_r$ has full Haar measure in $\Z_p$. \\

Let $\sigma$ be an element of $V$. As argued earlier, since $\sigma$ is non-zero, we write $\sigma$ as $p^{\mathrm{val}(\sigma)}u_\sigma$, for a unique element $u_\sigma$ in $\Z_p^\times$ and a non-negative integer $\mathrm{val}(\sigma)$. Consider the following commutative diagram:
\begin{align*}
\xymatrix{
\Z_p \ar[d]_{(\times u_\sigma^{-1})\circ \mathcal{T}_3 \circ (\times  u_\sigma)} \ar[r]^{\times u_\sigma}& \Z_p\ar[d]^{\mathcal{T}_3} \\
\Z_p & \Z_p \ar[l]^{\times u_\sigma^{-1}}
}	
\end{align*}
and the following set:
\begin{align}D_{u_\sigma} \coloneqq \left\{\alpha \in \Z_p, \ \text{ such that }  \left\{(u_\sigma^{-1}\circ \mathcal{T}_3 \circ u_\sigma)^{[n]} (\alpha)\right\}_{n=0}^{\infty} \text{ is equidistributed in }\Z_p  \right\}.
\end{align}
Notice that the map $(u_\sigma^{-1}\circ \mathcal{T}_3 \circ u_\sigma)^{[n]}$ equals $u_\sigma^{-1}\circ \mathcal{T}_3^{[n]} \circ u_\sigma$. As a result, we have
\begin{align}
D_{u_\sigma} =  \left\{\alpha \in \Z_p, \ \text{ such that }   \left\{u_\sigma^{-1}\circ \mathcal{T}_3^{[n]} \circ u_\sigma (\alpha) )\right\}_{n=0}^{\infty} \text{ is equidistributed in }\Z_p  \right\}.
\end{align}
Notice also that since $u_\sigma$ is a unit, the maps $\Z_p \xrightarrow {u_\sigma} \Z_p$ and $\Z_p \xrightarrow {u_\sigma^{-1}} \Z_p$ are continuous topological isomorphisms. As a result, for every continuous function $f:\Z_p \rightarrow \mathbb{R}$, the functions $f\circ u_\sigma:\Z_p \rightarrow \R$ and $f\circ u_\sigma^{-1}:\Z_p \rightarrow \R$ are also continuous functions. Since the maps $\Z_p \xrightarrow {u_\sigma} \Z_p$ and $\Z_p \xrightarrow {u_\sigma^{-1}} \Z_p$ also preserve the Haar measure, using the integration-by-substitution method, we  see that \[\int f \ \mathrm{d}\nu_{\mathrm{Haar}} = \int (f \circ u_\alpha) \ \mathrm{d}\nu_{\mathrm{Haar}} = \int (f\circ u_\alpha^{-1}) \  \mathrm{d}\nu_{\mathrm{Haar}}.\]
These observations let us conclude that a sequence $\{x_n\}$ is equidistributed in $\Z_p$ if and only if the sequence  $\{u_\sigma x_n\}$ is equidistributed in $\Z_p$ if and only if  the sequence $\{u_{\sigma}^{-1}x_n\}$ is equidistributed in $\Z_p$.
As a result, we can conclude that 
\begin{align}	
D_{u_\sigma} & = \left\{\alpha \in \Z_p, \ \text{ such that }    \left\{ \mathcal{T}_3^{[n]} (u_\sigma \alpha) )\right\}_{n=0}^{\infty}\text{ is equidistributed in }\Z_p  \right\}. 
\end{align}

 Observe that $\mathcal{T}_3^{[n]}(u_\sigma\alpha)=\mathcal{T}_3^{[n+ \mathrm{val}(\sigma)]}(\sigma\alpha)$.  Since equidistribution for a sequence is an asymptotic property, we can conclude that the sequence  $\{\mathcal{T}_3^{[n]}(\sigma \alpha)\}_{n=0}^\infty$   is equidistributed in $\Z_p$ if and only if $\left\{\mathcal{T}_3^{[n]}(u_{\sigma} \alpha)\right\}_{n=0}^\infty$  is equidistributed in $\Z_p$. As a result, we have 
\begin{align}	
D_{u_\sigma} & = \left\{\alpha \in \Z_p, \ \text{ such that }    \left\{ \mathcal{T}_3^{[n]} (\sigma \alpha) )\right\}_{n=0}^{\infty}\text{ is equidistributed in }\Z_p  \right\}. 
\end{align}

Using the implication (\ref{it:prop2}) $\iff$ (\ref{it:prop3}) of Proposition \ref{prop:nodependenceonrealcoordinates}, one has the following alternative description of $H_r$: 
\begin{align}
H_r =\left\{\alpha \in \Z_p, \text{ such that } \left(\alpha\beta_1,\cdots,\alpha\beta_r\right)\text{ is a generic point for } \mathcal{T}_3^r\right\}.	
\end{align}
Using Proposition \ref{prop:reductionstep}, we have that $\left(\alpha\beta_1,\cdots,\alpha\beta_r\right)$  is a generic point for  $\mathcal{T}_3^r$ if and only if for all $\sigma$ in $V$, the element $\sigma \alpha$ is a generic point for $\mathcal{T}_3$. That is, 
\begin{align}
H_r = \bigcap_{\sigma \in V} D_{u_\sigma}.
\end{align}

Since the maps $\Z_p \xrightarrow {u_\alpha} \Z_p$ and $\Z_p \xrightarrow {u_\alpha^{-1}} \Z_p$ preserve the Haar measure, one can directly check from the definition of an ergodic map, that since $\mathcal{T}_3:\Z_p \rightarrow \Z_p$ is an ergodic map, the map $u_\sigma^{-1}\circ \mathcal{T}_3 \circ u_\sigma:\Z_p \rightarrow \Z_p$ is also ergodic.

A simple application of Birkhoff's pointwise ergodic theorem given in \cite[Corollary 4.20]{MR2723325} lets us conclude that ergodicity of a self-map implies that almost every point is generic. Thus, using \cite[Corollary 4.20]{MR2723325}, we can conclude that $D_{u_\sigma}$ has full Haar measure (that is, has measure equal to $1$) in $\Z_p$. A countable intersection of sets with full Haar measure also has full Haar measure. As a result, the set $H_r$ also has full Haar measure in $\Z_p$. Proposition \ref{prop:equidistribution} follows.

\begin{remark}
Proposition \ref{prop:nodependenceonrealcoordinates}  involves establishing various equivalent characterisations of generic points for the self-maps of the $p$-adic solenoid and its fundamental domain. There are two main non-trivial characterisations. The first non-trivial characterisation of the proposition involves describing generic points solely from the dynamics on the $p$-adic side. This part of the proposition is described in Section \ref{subsec:31implication}. The main ingredient is Lemma \ref{lem:bernoullinodependence}  where a statement analogous to Proposition \ref{prop:nodependenceonrealcoordinates} is established for the symbolic space characterising generic points for the $2$-sided Bernoulli shift space in terms of generic points for the $1$-sided Bernoulli shift space. Note that the  $1$-sided Bernoulli shift space is identified with the $p$-adic side using $p$-adic expansions of $p$-adic integers. \\

On the other hand, equidistribution in a dynamical system encodes both topological and measure-theoretic information.  An equidistributed sequence in our setting, being countable, is a measure zero set. Although ultimately the symbolic space, the fundamental domain and the $p$-adic solenoid only differ by a measure zero set, we must be careful not to exclude any equidistributed sequences while characterising generic points involving measurably isomorphic systems.  The second non-trivial characterisation is described in section \ref{subsec:72implication}; the main point of that section is to characterise the generic points for the ergodic self-map on the fundamental domain in terms of the generic points for the ergodic self-map on the $p$-adic solenoid.
\end{remark}

\section{Compatibility of measures} \label{sec:compatibility}

Recall that the set $X_4$ equals $\{0,1,\cdots, p-1\}^{\Z}$. For each finite subset $I$ of $\Z$ and a map $b: I \rightarrow \{0,1,\cdots, p-1\}$, we can consider the set $I(b) = \{x \in X_4, \text{ such that } x_j = b(j)\}$. This set $I(b)$ is called a \textit{cylinder set}. The topology on $X_4$ is generated by these cylinder sets forming a basis of open sets, as $I$ varies over all finite subsets of $\Z$ and $b$ varies over all maps $I \rightarrow \{0,1,\cdots,p-1\}$. The sigma algebra on $X_4$ is the smallest sigma algebra containing all of the open sets generated by this topology. The measure $\nu_4$ on the Bernoulli scheme is given by the uniform probability vector $\left(1/p,1/p,\cdots, 1/p\right)$ on the sample space $\{0,1,\cdots,p-1\}$. As a result, the value of the measure $\nu_4$ on $I(b)$ is given by $\left({1}/{p}\right)^{|I|}$.

\begin{proposition} \label{prop:compatibility} 
Consider the following natural map:
\begin{align*}
\pi: X_4 & \twoheadrightarrow X_3,\\
\left( \cdots, a_{-n}, \cdots, a_{-2}, a_{-1}  \ \biggl\vert \  t_0, t_1, t_2, \cdots, t_n,\cdots\right) &\mapsto \left(\sum_{i=0}^\infty t_i p^i, \ \sum \limits_{j=1}^
\infty \dfrac{a_{-j}}{p^{j}} \right).
\end{align*} 
The following statements hold:
\begin{enumerate}
\item The map $\pi$ is a continuous surjective map.
\item The pushforward measure $\pi^*(\nu_4)$ on $X_3$ coincides with the measure $\nu_3$ on $X_3$.  
\end{enumerate}
\end{proposition}

\begin{proof}
Consider the following sets:
\[
Y_4 \coloneqq \left\{0,1,\cdots, p-1\right\}^{\Z_{\geq 0}}, \qquad 	Z_4 \coloneqq \left\{0,1,\cdots, p-1\right\}^{\Z_{< 0}}.
\] 
We shall assign topologies on $Y_4$ and $Z_4$ generated by the cylinder sets $I(b)$ forming the basis of open sets, as $I$ varies over finite subsets of $\Z_{\geq 0}$ (and $\Z_{<0}$ respectively) and $b$ varies over all maps $I \rightarrow \{0,1,\cdots, p-1\}$. One can also consider the smallest sigma algebras on $Y_4$ and $Z_4$ generated  by these topologies. As earlier, the corresponding measure $\nu_Y$ on $Y_4$ (and $\nu_Z$ on $Z_4$ respectively) assigns the value $\left({1}/{p}\right)^{|I|}$ to these cylinder sets. As earlier, this assignment corresponds to the uniform probability vector $\left(1/p,1/p,\cdots, 1/p\right)$ on the sample space $\{0,1,\cdots,p-1\}$. Note that as a topological space and as a measurable space, $X_4$ is naturally isomorphic to the product $Z_4 \times Y_4$. \\

Consider the following maps of topological and measurable spaces. 
 \begin{align*}
 	\pi_Z: \left(Z_4,\nu_Z\right) &\rightarrow \left([0,1],\nu_{\mathrm{Bor}}\right) \qquad  & \pi_Y: Y_4 &\rightarrow \left(\Z_p, \nu_{\mathrm{Haar}}\right)\\
 	\left(\cdots, a_{-n}, \cdots, a_{-2}, a_{-1}\right) &\mapsto \sum \limits_{j=1}^\infty \dfrac{a_{-j}}{p^{j}}, 
 	\qquad & \left(t_0, t_1, t_2, \cdots, t_n,\cdots\right) &\mapsto \sum_{i=0}^{\infty} t_i p^i. 
 \end{align*} 

To prove the proposition, it is enough to show that (i) both $\pi_Y$ and $\pi_Z$ are continuous surjective maps, and that (ii) the pushforward measures $\pi_Y^*(\nu_Y)$ and $\pi_Z^*(\nu_Z)$ coincide with $\nu_{\mathrm{Haar}}$  and $\nu_{\mathrm{Bor}}$ respectively. \\
	
For the fact that $\pi_Z$ is a continuous surjective map, see  \cite[Section 2.3, Theorem]{MR1760253}. Every element of $[0,1]$ has a unique $p$-ary expansion, except for the rational numbers whose denominators are powers of $p$; each such rational number has precisely two $p$-ary expansions (see \cite[Problem 44 in Chapter 1]{MR1013117}). For the fact that $\pi_Y$ is a continuous surjective map, see \cite[Section 2]{MR1760253}. In fact, $\pi_Y$ is an isomorphism of topological spaces. 	\\

	\begin{figure}
 \centering
\begin{tikzpicture}[thick]

  \node[circle,fill=black, inner sep=0pt,minimum size=1.5pt] (b) at (-2.2,0.25) {}; 
 \node[circle,fill=black, inner sep=0pt,minimum size=1.5pt] (b) at (-2.6,0.25) {};  
 \node[circle,fill=black, inner sep=0pt,minimum size=1.5pt] (b) at (-3,0.25) {}; 

 \filldraw[fill=white] (0.2,-0.1) rectangle (-0.5,0.6); 
 \filldraw[fill=white] (-0.6,-0.1) rectangle (-1.3,0.6);
 \filldraw[fill=white](-1.4,-0.1) rectangle (-2.1,0.6);

 \filldraw[fill=white] (-3.3,-0.1) rectangle (-4,0.6);
 
\node[circle,fill=black, inner sep=0pt,minimum size=5pt] (b) at (-4.4,0.25) {}; 
\node[circle,fill=black, inner sep=0pt,minimum size=5pt] (b) at (-4.8,0.25) {};  
 \node[circle,fill=black, inner sep=0pt,minimum size=5pt] (b) at (-5.2,0.25) {};  
 \node[circle,fill=black, inner sep=0pt,minimum size=5pt] (b) at (-5.6,0.25) {};  
 \node[circle,fill=black, inner sep=0pt,minimum size=5pt] (b) at (-6,0.25) {};  
 \node[circle,fill=black, inner sep=0pt,minimum size=5pt] (b) at (-6.4,0.25) {};  
 
 \node[] at (-0.13,0.2) {$a_{-1}$};
  \node[] at (-0.93,0.2) {$a_{-2}$};
  \node[] at (-1.72,0.2) {$a_{-3}$};
  \node[] at (-3.64,0.2) {$a_{-n}$};

   \draw [thick,domain=-45:45] plot ({cos(\x)}, {0.25+sin(\x)});
 \draw [thick,domain=135:226] plot ({-6+cos(\x)}, {0.25+ sin(\x)});

\end{tikzpicture}

 \vspace{2cm} 
\begin{tikzpicture}[
line/.style = {draw,thick, 
               shorten >=-2pt, shorten <=-2pt}
                     ]
\draw (0,0) -- (10,0);
\foreach \i in {0,1} 
\draw (0,0.15) -- ++ (0,-0.3) node[below] {$0$}; 
\draw (10,0.15) -- ++ (0,-0.3) node[below] {$1$}; 
\draw[line, ultra thick, {Circle[length=7pt,fill=white]}-{Circle[length=7pt, fill=white]}]  (3,0) -- (5,0);
\draw (3,-0.5)  node[below] {$\dfrac{a}{p^{n}}$}; 
\draw (5,-0.5)  node[below] {$\dfrac{a+1}{p^{n}}$}; 
  \end{tikzpicture}	
  \caption{Cylinder sets and open intervals}
\end{figure}

To show that the pushforward  measure $\pi_Z^*(\nu_Z)$ coincides with $\nu_{\mathrm{Bor}}$,  consider the following collection of open intervals in $[0,1]$:
\[\mathfrak{C}_1 \coloneqq \left\{\left(\frac{a}{p^n}, \frac{a+1}{p^n}\right), \text{for all } n \in \Z_{\geq 0} \text{ and } 0 \leq a < p^n, \text{ with } \mathrm{gcd}(a,p) =1 \right\},	\]
along with the following collection of singleton sets of $[0,1]$:
\[\mathfrak{C}_2 \coloneqq \left\{\left\{\frac{a}{p^n}\right\}, \text{for all } n \in \Z_{\geq 0} \text{ and } 0 \leq a \leq p^n \right\}.	\]

We let $\mathfrak{C}$ equal $\mathfrak{C}_1  \cup  \mathfrak{C}_2 \cup \{\emptyset\}$. We will show that $\mathfrak{C}$ is a separating class (in the sense of Billingsley's book \cite{MR1700749}) and hence, it is enough to check that $\pi_Z^*(\nu_Z)$ coincides with $\nu_{\mathrm{Bor}}$ for all the elements of $\mathfrak{C}$. Although the arguments are fairly elementary, we provide some details for the reader's convenience. \\

Firstly, we show $\mathfrak{C}$  is a $\pi$-system (in the sense of Billingsley's book \cite{MR2893652}). That is, we claim that $\mathfrak{C}$ is closed under finite intersections. To see this, since  $\mathfrak{C}_2$ is a collection of singleton sets, it is enough to show that the intersection of two open intervals of the form $\left(\frac{a}{p^n}, \frac{a+1}{p^n}\right)$ and $\left(\frac{b}{p^s}, \frac{b+1}{p^s}\right)$ also belongs to $\mathfrak{C}$. 
\begin{itemize}
\item If $\frac{a}{p^n} = \frac{b}{p^s}$, then $a=b$ and $p^n=p^s$; these two intervals must be equal to each other and hence their intersection belongs to $\mathfrak{C}$.	
\item  Without loss of generality, we may assume that $\frac{a}{p^n} < \frac{b}{p^s}$. In this case, if $\frac{a+1}{p^n} \leq \frac{b}{p^s}$, then the intersection of the two intervals would equal the empty set, which also belongs to $\mathfrak{C}$. 
\item The last case to consider is $\frac{a}{p^n} < \frac{b}{p^s} < \frac{a+1}{p^n}$. In this last case, we cannot have $n \geq s$, since otherwise we would get $a < p^{n-s}b < a+1$. However, there are no integers strictly between $a$ and $a+1$.  Therefore, in this last case, we must have $0 \leq n < s$. As a result, we must then have $p^{s-n}a < b < p^{s-n}(a+1)$. Since $s-n$ is at least one, the difference $p^{s-n}(a+1) - p^{s-n}a$ must be at least $p$. Consequently, we must have $p^{s-n}a < b+1 < p^{s-n}(a+1)$. In other words, in this last case, we have 
\[\frac{a}{p^n} < \frac{b}{p^s} < \frac{b+1}{p^s} < \frac{a+1}{p^n}.\]
The intersection of the two intervals $\left(\frac{a}{p^n}, \frac{a+1}{p^n}\right)$ and $\left(\frac{b}{p^s}, \frac{b+1}{p^s}\right)$ must equal $\left(\frac{b}{p^s}, \frac{b+1}{p^s}\right)$ and hence also belong to $\mathfrak{C}$.
\end{itemize}
 
Secondly, we claim that the sigma algebra $\Sigma_\mathfrak{C}$ generated by the elements of $\mathfrak{C}$ coincides with the Borel sigma algebra of $[0,1]$. Since $\mathfrak{C}$ is a collection of open intervals and closed (singleton) sets in $[0,1]$, the Borel sigma algebra of $[0,1]$ contains $\Sigma_\mathfrak{C}$. To show that reverse inclusion, we must show that the sigma algebra $\Sigma_\mathfrak{C}$ contains all the open sets of $[0,1]$. It will be sufficient to show that every open set can be written as a union of sets of $\mathfrak{C}$, since $\mathfrak{C}$ itself is a countable collection of sets of $[0,1]$.  Let $U$ be an open set  in $[0,1]$. We must show that for every $x$ in $[0,1]$, there exists a set $V$ in $\mathfrak{C}$ such that $x \in V \subset U$.
\begin{itemize}
\item Consider a point $x$ in $[0,1]$. If $x$ is of the form $\frac{a}{p^n}$ for some $n \in \Z_{\geq 0}$ and $0 \leq a \leq p^n$, then we may choose $V$ to be the set $\{x\}$, since $\{x\}$ belongs to $\mathfrak{C}_2$.
\item Consider a point $x$ in $[0,1]$ which is not of the form  $\frac{a}{p^n}$  as above. The $p$-ary expansion of $x$ is then unique and the $p$-ary expansion does not end with the repeating sequence $0,0,0,\cdots,$ or the sequence $p-1,p-1,p-1,\cdots,$. Furthermore, since $U$ is open, there exists an integer $n$ such that (i) $(x-\frac{1}{p^n},x+\frac{1}{p^n}) \subset U$ and (ii) the $n$-th digit in the $p$-ary expansion of $x$ is not zero. For such an $n$, we let $\frac{a}{p^n}$ denote the truncation of the $p$-ary expansion of $x$ at the $n$-th digit. Observe that since the $n$-th digit is not $0$, the numerator $a$ is co-prime to $p$. Since the $p$-ary expansion does not end with the repeating sequence $0,0,0,\cdots,$, we have $x > \frac{a}{p^n}$. Since the $p$-ary expansion does not end with the repeating sequence $p-1,p-1,p-1,\cdots,$, we have $x < \frac{a+1}{p^n}$. Combining these observations, we have 
\[x-\frac{1}{p^n} < \frac{a}{p^n} < \frac{a+1}{p^n} < x+\frac{1}{p^n}.\]
In other words, $x \in \left(\frac{a}{p^n}, \frac{a+1}{p^n}\right) \subset (x-\frac{1}{p^n},x+\frac{1}{p^n}) \subset U$. In this case, we may take $V$ to be $\left(\frac{a}{p^n}, \frac{a+1}{p^n}\right)$.
\end{itemize}

 By combining these two claims and applying \cite[Lemma 1.9.4]{MR2267655}, to show that the pushforward  measure $\pi_Z^*(\nu_Z)$ coincides with $\nu_{\mathrm{Bor}}$, it is enough to check that these probability measures agree on the sets in $
\mathfrak{C}$.  It is not hard to check that singleton sets in the Bernoulli scheme $Y_4$, for the uniform probability vector $\left(1/p,1/p,\cdots, 1/p\right)$, are closed and have measure zero. As a result, the two measures $\pi_Z^*(\nu_Z)$ and $\nu_{\mathrm{Bor}}$ agree on all the sets in $\mathfrak{C}_2$ (and are equal to zero).  

Consider an interval $\left(\frac{a}{p^n}, \frac{a+1}{p^n}\right)$ in $\mathfrak{C}_1$. Consider also the $p$-ary expansion of $\frac{a}{p^n}$:
\[\frac{a}{p^n} = \frac{a_{-1}}{p^{1}} + \cdots + \frac{a_{-n}}{p^n}.\]  \\
Let $I$ be the finite subset $\{-1,-2,\cdots, -n\}$ of $\Z_{<0}$. Let $b$ be the assignment that sends each $i$ in $I$ to $a_{i}$ in $\{0,1,\cdots,p-1\}$. Let $z_1, z_2$ be the elements $\left(,\cdots, 0, 0, 0, 0, a_{-n},a_{-n+1}, \cdots, a_{-1} \right)$ and $\left(,\cdots, p-1, p-1, p-1, p-1, a_{-n},a_{-n+1}, \cdots, a_{-1} \right)$ respectively.  The preimage of $\left(\frac{a}{p^n}, \frac{a+1}{p^n}\right)$ under $\pi_Z$ equals 
\[I(b) - \{z_1,z_2\}.
\]
Combining these observations, we can directly verify that the values of the measures $\pi_Z^*(\nu_Z)$ and $\nu_{\mathrm{Bor}}$ coincide on  $\left(\frac{a}{p^n}, \frac{a+1}{p^n}\right)$ and are both equal to $1/p^n$. Therefore, the measures $\pi_Z^*(\nu_Z)$ and $\nu_{\mathrm{Bor}}$ agree. \\

The computation to show that the pushforward measure $\pi_Y^*(\nu_Y)$ coincides with $\nu_{\mathrm{Haar}}$ is almost exactly similar to the discussion above. We skip it for the sake of brevity. It would be enough to that these two measures coincide for the open metric balls in $\Z_p$, since (i) these open metric balls generate the sigma algebra for the Haar measure on $\Z_p$ and (ii) these open metric balls form a $\pi$-system because, in the non-archimedean metric, two metric balls are either disjoint or one of them is contained in the other. 
\end{proof}

\section{Transferring dynamics to the $p$-adic side: Proof of Proposition \ref{prop:nodependenceonrealcoordinates}} \label{sec:proofpropnodependence}

The implications (\ref{it:prop1}) $\implies$ (\ref{it:prop2}), (\ref{it:prop4}) $\implies$ (\ref{it:prop6}) and (\ref{it:prop5}) $\implies$ (\ref{it:prop7}) of Proposition \ref{prop:nodependenceonrealcoordinates} is clear. \\

Just as in Section \ref{subsec:proofimplication}, the implication (\ref{it:prop2}) $\implies$ (\ref{it:prop3}) of Proposition \ref{prop:nodependenceonrealcoordinates} follows by observing that the measure $\nu_{\mathrm{Haar}}^r$ on $\Z_p^r$ coincides with the pushforward measure of $(\nu_{\mathrm{Haar}} \times \nu_{\mathrm{Bor}})^r$ under the natural continuous projection map $(\Z_p\times[0,1])^r \rightarrow \Z_p^r$ and then applying \cite[Theorem 3.6.1]{MR2267655}.  The implication (\ref{it:prop1}) $\implies$ (\ref{it:prop4}) follows by a similar reasoning by considering the natural map $X_3^r \rightarrow X_2^r$. The implication (\ref{it:prop4}) $\iff$ (\ref{it:prop5}) and (\ref{it:prop6}) $\iff$ (\ref{it:prop7}) are also straightforward because the natural topological (and measurable) isomorphism $X_2^r \cong X_1^r$ commutes with the ergodic maps $T_2^r$ and $T_1^r$. To complete the proof of Proposition \ref{prop:nodependenceonrealcoordinates}, see the proof of implications given in sections \ref{subsec:31implication} and \ref{subsec:72implication}.

\subsection{The implication (\ref{it:prop3}) $\implies$ (\ref{it:prop1})}\label{subsec:31implication}
To see the implication (\ref{it:prop3}) $\implies$ (\ref{it:prop1}) of Proposition \ref{prop:nodependenceonrealcoordinates}, consider the following commutative diagram of measurable spaces:
\begin{align*}
\xymatrix{
X_4^r \coloneqq \left(\left\{0,1,\cdots, p-1\right\}^{\Z}\right)^r \ar[d] \ar@{->>}[r]&  \left((\Z_p \times [0,1])^r, (\nu_{\mathrm{Haar}} \times \nu_{\mathrm{Bor}})^r \right) \ar[d] \\
 Y_4^r \coloneqq \left(\left\{0,1,\cdots, p-1\right\}^{\Z_{\geq 0}}\right)^r \ar@{->>}[r] &  \left(\Z_p^r, \nu_{\mathrm{Haar}}^r\right)
 }		
\end{align*}
Recall that on the left hand  side of the commutative diagram, the two-sided and the one-sided  Bernoulli schemes
\[\left\{0,1,\cdots, p-1\right\}^{\Z} \stackrel{T_4}{\longrightarrow}\left\{0,1,\cdots, p-1\right\}^{\Z}, \qquad \left\{0,1,\cdots, p-1\right\}^{\Z_{\geq 0}} \stackrel{\mathcal{T}_4}{\longrightarrow} \left\{0,1,\cdots, p-1\right\}^{\Z_{\geq 0}},   \]
  are both equipped with the probability measures induced by the uniform probability vector $\left(1/p,1/p,\cdots, 1/p\right)$ on the sample space $\{0,1,\cdots,p-1\}$. 
  
  As observed in the proof of Proposition \ref{prop:compatibility}, the bottom horizontal map $Y_4^r \stackrel{\cong}{\longrightarrow}  \Z_p^r$ is an isomorphism of topological and measurable spaces. The maps $T_4^r$ and $\mathcal{T}_4^r$, on $X_4^r$ and $Y_4^r$ respectively, along with the maps $T_3^r$ and $\mathcal{T}_3^r$, on $(\Z_p\times[0,1])^r$ and $\Z_p^r$ respectively, are compatible with the maps in the commutative diagram given above.  Consequently, we make the following observations: 
   \begin{itemize} 
  \item $(y_1,\cdots,y_r)$ is a generic point for $\mathcal{T}_4^r$ if and only if its image, under the bottom horizontal map, is a generic point for $\mathcal{T}_3^r$.
  \item Using Proposition \ref{prop:compatibility}, we can conclude that the pushforward of the measure $\nu_4^r$ coincides with $(\nu_{\mathrm{Haar}}\times \nu_{\mathrm{Bor}})^r$.  By applying \cite[Theorem 3.6.1]{MR2267655}, we can conclude that if the element $\big(\left(z_1,y_1\right),\left(z_2,y_2\right),\cdots,\left(z_r,y_r\right)\big)$ is a generic point for $T_4^r$, then its image, under the top horizontal map of the commutative diagram, is a generic point for $T_3^r$.
  \end{itemize}
 As in Proposition \ref{prop:compatibility}, we let $Z_4$ denote  $\left\{0,1,\cdots, p-1\right\}^{\Z_{< 0}}$. One can use these observations to conclude that, to prove the implication (\ref{it:prop3}) $\implies$ (\ref{it:prop1}) of Proposition \ref{prop:nodependenceonrealcoordinates}, it suffices to work with Bernoulli schemes and use the implication of (\ref{it:3}) $\implies$ (\ref{it:1}) of Lemma \ref{lem:bernoullinodependence}. The key observation is that we are working with a \textit{left} Bernoulli shift.

\begin{lemma}\label{lem:bernoullinodependence}
	The following statements are equivalent:
\begin{enumerate}
\item\label{it:1} For every $(z_1,\cdots,z_r)$ in $Z_4^r$, the element $\big(\left(z_1,y_1\right),\left(z_2,y_2\right),\cdots,\left(z_r,y_r\right)\big)$ is a generic point for $T_4^r$.
\item\label{it:2} There exists a $(z_1,\cdots,z_r)$ in $Z_4^r$ such that the element $\big(\left(z_1,y_1\right),\left(z_2,y_2\right),\cdots,\left(z_r,y_r\right)\big)$ is a generic point for $T_4^r$.
\item\label{it:3}  $(y_1,\cdots,y_r)$ is a generic point for $\mathcal{T}_4^r$.
\end{enumerate}
\end{lemma}

\begin{proof}
The implication (\ref{it:1}) $\implies$ (\ref{it:2}) is clear. \\

To prove the implications  (\ref{it:2}) $\implies$ (\ref{it:3}) $\implies$ (\ref{it:1}), we will first show that the finite product of cylinder sets for the Bernoulli schemes, along with the empty set, form a separating class (in the sense of Billingsley's book \cite{MR1700749}).
Observe that the collection of (a finite product of) cylinder sets is countable. To see this, one simply needs to use the fact that the collection of finite subsets of a countable set is countable. Since the topology and the sigma algebra on the Bernoulli schemes are generated by the product of cylinder sets, it is enough to show that the cylinder sets form a $\pi$-system (in the sense of Billingsley's book \cite{MR2893652}). It is enough to consider the case $r=1$. To see that the intersection of two cylinder sets is again a cylinder set, we let $I_1(b_1)$ and $I_2(b_2)$ denote two cylinder sets. If there exists an $i \in I_1 \cap I_2$, such that $b_1(i) \neq b_2(i)$, then the intersection of the  $I_1(b_1)$ and $I_2(b_2)$ is the empty set. Otherwise, the intersection of $I_1(b_1)$ and $I_2(b_2)$ equals the cylinder set $J(c)$, where $J$ is the union of the finite sets $I_1$ and $I_2$, and the function $c:J \rightarrow \{0,1,\cdots, p-1\}$ is defined to be $b_1(i)$ if $i \in I_1$ and $b_2(i)$ if $i \in I_2$. 

Observe that the inverse image of product of cylinder sets under the map $X_4^r \rightarrow Y_4^r$ is again a product of cylinder sets. Combining the observations in the previous paragraph and applying \cite[Lemma 1.9.4]{MR2267655}, one can conclude that the pushforward of the Bernoulli measure $\nu_4^r$ coincides with the Bernoulli measure on $Y_4^r$.  

Observe that the product of cylinder sets is given by the inverse image of a continuous map from $X_4^r$ (and $Y_4^r$ respectively) to a finite discrete set. Consequently, the product of the cylinder sets must be both open and closed in $X_4^r$ (and $Y_4^r$ respectively). As a result, the boundary of a product of the cylinder sets is the empty set. The product of the cylinder sets is thus a \textit{$P$-continuity set}, in the sense of Billingsley's book \cite[Section 2]{MR1700749}. \\

To prove the implication (\ref{it:2}) $\implies$ (\ref{it:3}), assume that there exists a $(z_1,\cdots,z_r)$ in $Z_4^r$ such that the element $\big(\left(z_1,y_1\right),\left(z_2,y_2\right),\cdots,\left(z_r,y_r\right)\big)$ is a generic point for $T_4^r$. In other words, the set $\left\{(T_4^r)^n\left(\big(\left(z_1,y_1\right),\left(z_2,y_2\right),\cdots,\left(z_r,y_r\right)\big)\right)\right\}_{n}$ is equidistributed in $X_4^r$. Let  $\varpi$ denote the natural projection map $X_4^r \rightarrow Y_4^r$. Just as in Section \ref{subsec:proofimplication}, applying \cite[Theorem 3.6.1]{MR2267655} lets us conclude that the set $\left\{(\varpi \circ T_4^r)^n\left(\big(\left(z_1,y_1\right),\left(z_2,y_2\right),\cdots,\left(z_r,y_r\right)\big)\right)\right\}_{n}$ is equidistributed in $Y_4$.
Using the following commutative diagram
\begin{align*}
\xymatrix{
X_4^r   \ar[d]_{\varpi} \ar[r]^{T_4^r}&  X_4^r  \ar[d]^{\varpi} \\
 Y_4^r \ar[r]^{\mathcal{T}_4^r} &  Y_4^r
 }		
\end{align*}
we have  
\[(\varpi \circ T_4^r)^{[n]}\left(\big(\left(z_1,y_1\right),\left(z_2,y_2\right),\cdots,\left(z_r,y_r\right)\big)\right) = (\mathcal{T}_4^r)^{[n]}\left((y_1,y_2,\cdots, y_r)\right).\]
That is, the set $\{(\mathcal{T}_4^r)^{[n]}\left((y_1,y_2,\cdots, y_r)\right)\}_n$ is equidistributed in $Y_4$. In other words, $(y_1,\cdots,y_r)$ is a generic point for $\mathcal{T}_4^r$. This proves the implication (\ref{it:2}) $\implies$ (\ref{it:3}). \\

To prove the implication (\ref{it:3}) $\implies$ (\ref{it:1}), consider the following equation:
\begin{align} \label{eq:equidef}
\lim_{M \rightarrow \infty} \dfrac{1}{M}\sum_{n=1}^M f\left(x_n\right) = 	\int_X f(x) d \nu_X
\end{align}
Here, $X \in \{X_4^r, \ Y_4^r\}$ and $\{x_n\}$ is a sequence of points in $X$. The following statements are equivalent:
\begin{enumerate}[label=(\roman*)]
\item\label{it:roman1} Equation (\ref{eq:equidef}) holds for all continuous functions $f: X \rightarrow [0,1]$.
\item\label{it:roman2} Equation (\ref{eq:equidef}) holds for the indicator function $\mathds{1}_{I_1(b_1) \times \cdots \times I_r(b_r)}: X \rightarrow [0,1]$, for all cylinder sets $I_1(b_1), \cdots, I_r(b_r)$ in $X$,  with $b_i: I_i \rightarrow \{ 0, \cdots, p-1 \}$.
\end{enumerate}

The implication \ref{it:roman1} $\implies$ \ref{it:roman2} follows by an application of the Portmanteau theorem (\cite[Theorem 2.1]{MR1700749}) since the product of cylinders sets is a $P$-continuity set. 	The reverse implication \ref{it:roman2} $\implies$ \ref{it:roman1} follows by (\cite[Theorem 2.2]{MR1700749}) since  products of cylinder sets form a $\pi$-system class and since there are only a countable collection of such products.  \\

Suppose $(y_1,\cdots,y_r)$ is a generic point for $\mathcal{T}_4^r$. That is, for all cylinder sets $J_1(c_1), \cdots, J_r(c_r)$, we have 
\begin{align}\label{eq:equidistributed1side}
	\lim_{M \rightarrow \infty} \dfrac{1}{M}\sum_{m=1}^M \mathds{1}_{J_1(c_1) \times \cdots \times J_r(c_r)}\left((\mathcal{T}_4^r)^{[m]}(y_1,\cdots,y_r)\right) = \prod_{i=1}^r \left(\dfrac{1}{p}\right)^{|J_i|}.
\end{align}
We need to show that for every $(z_1,\cdots,z_r)$ in $Z_4^r$, the element $\big(\left(z_1,y_1\right),\left(z_2,y_2\right),\cdots,\left(z_r,y_r\right)\big)$ is a generic point for $T_4^r$. That is, for all cylinder sets in $X_4$, we need to show
\begin{align}\label{eq:toshoweq}
	\lim_{M \rightarrow \infty} \dfrac{1}{M}\sum_{n=1}^M \mathds{1}_{I_1(b_1) \times \cdots \times I_r(b_r)}\left((T_4^r)^{[n]}\big(\left(z_1,y_1\right),\left(z_2,y_2\right),\cdots,\left(z_r,y_r\right)\big)\right) \stackrel{?}{=} \prod_{i=1}^r \left(\dfrac{1}{p}\right)^{|I_i|}.
\end{align}

 \vspace{0.5cm}
	\begin{figure}
 \centering
\begin{tikzpicture}[thick]

\node[circle,fill=black, inner sep=0pt,minimum size=5pt] (b) at (2.4,0.25) {};
\node[circle,fill=black, inner sep=0pt,minimum size=5pt] (b) at (2.8,0.25) {};
\node[circle,fill=black, inner sep=0pt,minimum size=5pt] (b) at (3.2,0.25) {};
\node[circle,fill=black, inner sep=0pt,minimum size=5pt] (b) at (3.6,0.25) {};
\node[circle,fill=black, inner sep=0pt,minimum size=5pt] (b) at (4,0.25) {};
\node[circle,fill=black, inner sep=0pt,minimum size=5pt] (b) at (4.4,0.25) {};
\node[circle,fill=black, inner sep=0pt,minimum size=5pt] (b) at (4.8,0.25) {};
\node[circle,fill=black, inner sep=0pt,minimum size=5pt] (b) at (5.2,0.25) {};
\node[circle,fill=black, inner sep=0pt,minimum size=5pt] (b) at (5.6,0.25) {};
\fill[gray] (0.5,0) rectangle (0,0.5);
\fill[gray] (0.7,0) rectangle (1.2,0.5);
\fill[gray] (1.4,0) rectangle (1.9,0.5);

\node[circle,fill=black, inner sep=0pt,minimum size=5pt] (b) at (-1.5,0.25) {};
\node[circle,fill=black, inner sep=0pt,minimum size=5pt] (b) at (-1.1,0.25) {};
\node[circle,fill=black, inner sep=0pt,minimum size=5pt] (b) at (-0.7,0.25) {};
 \fill[palesilver] (-2,0) rectangle (-2.5,0.5);
\fill[palesilver] (-2.7,0) rectangle (-3.2,0.5);
\fill[palesilver] (-3.4,0) rectangle (-3.9,0.5);
\node[circle,fill=black, inner sep=0pt,minimum size=5pt] (b) at (-4.4,0.25) {}; 
\node[circle,fill=black, inner sep=0pt,minimum size=5pt] (b) at (-4.8,0.25) {};  
 \node[circle,fill=black, inner sep=0pt,minimum size=5pt] (b) at (-5.2,0.25) {};  
 \node[circle,fill=black, inner sep=0pt,minimum size=5pt] (b) at (-5.6,0.25) {};  
 \node[circle,fill=black, inner sep=0pt,minimum size=5pt] (b) at (-6,0.25) {};  
 \node[circle,fill=black, inner sep=0pt,minimum size=5pt] (b) at (-6.4,0.25) {};  
 
   \node at (1,-0.8) {$J(c)$};
   \node at (-3,-0.8) {$I(b)$};
   \draw[ultra thick,<-] (-2,1.5) -- (0,1.5); 
    \node at (-1,2) {Bernoulli Shift};

\draw[|-|] (-3.9,-1.5) -- (-0.3,-1.5); 
  \node at (-2,-2) {$\alpha$};
 
   \draw [thick,domain=-45:45] plot ({5+cos(\x)}, {0.25+sin(\x)});
 \draw [thick,domain=135:226] plot ({-6+cos(\x)}, {0.25+ sin(\x)});

\draw[-] (-0.3,1) -- (-0.3,-1);

\draw [decorate, 
    decoration = {brace,raise=5pt}] (-2,0) --  (-3.9,0);
 
 \draw [ultra thick, decorate, 
    decoration = {brace,raise=5pt}] (5.5,-2.2) --  (0,-2.2); 
      \node at (3,-3) {$Y_4$};
 
  \draw [ultra thick,decorate, 
    decoration = {brace,raise=5pt}] (5.5,-3.5) --  (-6,-3.5); 
      \node at (-0.5,-4.2) {$X_4$};
 
\draw [decorate, 
    decoration = {brace,raise=5pt}] (1.9,0) --  (0,0);
 
\end{tikzpicture}
\caption{Visualizing the cylinder sets}
\end{figure}
 \vspace{0.5cm}

Let $\alpha \coloneqq \min \{x \in I_i, \text{ for } 1 \leq i \leq r\}$. For each $1 \leq i \leq r$, we define $J_i$ to be the set $\{x + \alpha, \text{ for all } x \in I_i\}$.  We define $c_i$ to be the function $c_i(y)= b_i(y-\alpha)$. Note that 
\begin{align} \label{eq:equalityofsetsizes}
|I_i| = |J_i|	
\end{align}
Since we are working with left Bernoulli shifts, for each point $\delta$ in $X_4^r$, we have 
\begin{align}\label{eqq:temp1}
	 \mathds{1}_{I_1(b_1) \times \cdots \times I_r(b_r)} \bigg((T_4^r)^{[\alpha]}(\delta)\bigg)  = \mathds{1}_{J_1(c_1) \times \cdots \times J_r(c_r)}\bigg(\delta\bigg).
\end{align}\label{eqq:temp2}
By our choice of $\alpha$, each of the sets $J_i$ is a finite subset of $\Z_{\geq 0}$.  We may view each of these cylinder sets in $Y_4^r$ view the natural projection map $X_4^r \rightarrow Y_4^r$. As a result, we have 
\begin{align}
\mathds{1}_{J_1(c_1) \times \cdots \times J_r(c_r)}\bigg(\delta\bigg)  = \mathds{1}_{J_1(c_1) \times \cdots \times J_r(c_r)}\bigg(\varpi\circ\delta\bigg).
\end{align}

By specializing $\delta$ to be of the form $(T_4^r)^{[m]}\big(\left(z_1,y_1\right),\left(z_2,y_2\right),\cdots,\left(z_r,y_r\right)\big)$ and combining equations (\ref{eqq:temp1}) and (\ref{eqq:temp2}), we obtain the following equality:
\begin{align}\label{eqq:temp3}
\notag
  & \mathds{1}_{I_1(b_1) \times \cdots \times I_r(b_r)} \bigg((T_4^r)^{[m+\alpha]} \big(\left(z_1,y_1\right),\left(z_2,y_2\right),\cdots,\left(z_r,y_r\right)\big)\bigg)  \\ &  = \mathds{1}_{J_1(c_1) \times \cdots \times J_r(c_r)} \bigg((\mathcal{T}_4^r)^{[m]}\big(y_1,\cdots,y_r\big)\bigg). \end{align}

Since equidistribution is an asymptotic property, one can discard finitely many terms in a sequence. Establishing equation (\ref{eq:toshoweq}) is equivalent to establishing the following equality:
\begin{align}\label{eq:toshoweq2}
	\lim_{M \rightarrow \infty} \dfrac{1}{M}\sum_{n=\alpha+1}^M \mathds{1}_{I_1(b_1) \times \cdots \times I_r(b_r)}\bigg((T_4^r)^{[n]}\big(\left(z_1,y_1\right),\left(z_2,y_2\right),\cdots,\left(z_r,y_r\right)\big)\bigg) \stackrel{?}{=} \prod_{i=1}^M \left(\dfrac{1}{p}\right)^{|I_i|}.
\end{align}
By letting $m=n-\alpha$, establishing equation (\ref{eq:toshoweq2}) is equivalent to establishing the following equality: 
\begin{align}\label{eq:toshoweq3}
	\lim_{M \rightarrow \infty} \dfrac{1}{M}\sum_{m=1}^{M-\alpha} \mathds{1}_{I_1(b_1) \times \cdots \times I_r(b_r)}\left((T_4^r)^{[m+\alpha]}\big(\left(z_1,y_1\right),\left(z_2,y_2\right),\cdots,\left(z_r,y_r\right)\big)\right) \stackrel{?}{=} \prod_{i=1}^M \left(\dfrac{1}{p}\right)^{|I_i|}.
\end{align}
Combining equations (\ref{eq:equalityofsetsizes}) and (\ref{eqq:temp3}), we can assert that establishing equation (\ref{eq:toshoweq3}) is equivalent to establishing the following equality:  
\begin{align}\label{eq:toshoweq4}
	\lim_{M \rightarrow \infty} \dfrac{1}{M}\sum_{m=1}^{M-\alpha}  \mathds{1}_{J_1(c_1) \times \cdots \times J_r(c_r)} \left((\mathcal{T}_4^r)^{[m]}\big(y_1,y_2,\cdots,y_r\big)\right) \stackrel{?}{=} \prod_{i=1}^M \left(\dfrac{1}{p}\right)^{|J_i|}.
\end{align}

Equation (\ref{eq:toshoweq4}) follows from equation (\ref{eq:equidistributed1side}) since  $(y_1,\cdots,y_r)$ is a generic point for $\mathcal{T}_4^r$. This concludes the proof of the implication (\ref{it:3}) $\implies$ (\ref{it:1}) and hence the proof of Lemma \ref{lem:bernoullinodependence}. 
This also concludes the proof of the implication (\ref{it:prop3}) $\implies$ (\ref{it:prop1}) in Proposition \ref{prop:nodependenceonrealcoordinates}.
 \end{proof}

\subsection{The implication (\ref{it:prop7}) $\implies$ (\ref{it:prop2})}\label{subsec:72implication}

We begin with a general setup. Let $\mathcal{X}$, $\mathcal{Y}$ be two compact topological spaces. Consider two dynamical systems $(\mathcal{X},m_\mathcal{X},\mathfrak{T}_\mathcal{X})$ and $(\mathcal{Y},m_\mathcal{Y},\mathfrak{T}_\mathcal{Y})$  along with a continuous surjection $\varpi: \mathcal{X} \rightarrow \mathcal{Y}$ such that the following diagram is commutative: 
\begin{align}\label{commdiag:generalsetup}
\xymatrix{
\mathcal{X}\ar[d]^{\varpi}  \ar[r]^{\mathfrak{T}_\mathcal{X}}	& \mathcal{X} \ar[d]^{\varpi} \\
\mathcal{Y} \ar[r]^{\mathfrak{T}_\mathcal{Y}}& 	\mathcal{Y} 
}
\end{align}
We assume that this setup satisfies the following additional conditions:
\begin{enumerate}[label=(\roman*), ref=\roman*]
\item\label{item:gen1} The  measures $m_\mathcal{X}$ and $m_\mathcal{Y}$ are probability measures defined on the Borel sigma algebras of $\mathcal{X}$ and $\mathcal{Y}$ respectively. 
\item\label{item:gen2}  The measure $m_\mathcal{Y}$ equals the pushfoward measure $\varpi^*(m_\mathcal{X})$.
\item\label{item:gen3} There exists sets $\mathcal{X}_0$, $\mathcal{Y}_0$  inside $\mathcal{X}$, $\mathcal{Y}$ respectively satisfying the following hypotheses: 
\begin{enumerate}
\item \label{item:gen3a} $\mathcal{X}_0$, $\mathcal{Y}_0$ are closed in $\mathcal{X}$, $\mathcal{Y}$ respectively. Both of them have measure zero.
\item\label{item:gen3b} The map $\varpi: \mathcal{X}_0^{\mathsf{c}}  \rightarrow \mathcal{Y}_0^{\mathsf{c}} $ is a homeomorphism. We let $\varphi$ denote the continuous inverse. 
\item\label{item:gen3c} The transformations $\mathcal{T}_\mathcal{X}$, $\mathcal{T}_\mathcal{Y}$ preserve $\mathcal{X}_0^{\mathsf{c}}$ and $\mathcal{Y}_0^{\mathsf{c}}$ respectively. That is, we have maps \[\mathcal{T}_\mathcal{X}: \mathcal{X}_0^{\mathsf{c}} \rightarrow \mathcal{X}_0^{\mathsf{c}}, \quad \mathcal{T}_\mathcal{Y}: \mathcal{Y}_0^{\mathsf{c}} \rightarrow  \mathcal{Y}_0^{\mathsf{c}}.\] \end{enumerate}
\end{enumerate}

 Here, $\mathcal{X}_0^{\mathsf{c}}$, $\mathcal{Y}_0^{\mathsf{c}}$ denote $\mathcal{X} \setminus \mathcal{X}_0$ and $\mathcal{Y} \setminus \mathcal{Y}_0$ respectively. It is not hard to see using  the conditions (\ref{item:gen3b}),(\ref{item:gen3c}), that the commutative diagram (\ref{commdiag:generalsetup}) can be extended to the following commutative diagram:
\begin{align}\label{commdiag:generalsetup2}
\xymatrix{
\mathcal{X}_0^{\mathsf{c}}\ar[d]^{\varpi}  \ar[r]^{\mathfrak{T}_\mathcal{X}}	& \mathcal{X}_0^{\mathsf{c}}\ar[d]^{\varpi} \\
\mathcal{Y}_0^{\mathsf{c}} \ar[d]^{\varphi}  \ar[r]^{\mathfrak{T}_\mathcal{Y}}& 	\mathcal{Y}_0^{\mathsf{c}} \ar[d]^{\varphi}  \\
\mathcal{X}_0^{\mathsf{c}} \ar[r]^{\mathfrak{T}_\mathcal{X}}& 	 \mathcal{X}_0^{\mathsf{c}}
}
\end{align}

\begin{lemma}
Suppose that $A$ is a $P$-continuity set in $\mathcal{X}$. Then, $\varpi\big(A \cap \mathcal{X}_0^{\mathsf{c}}\big)$ is a $P$-continuity set in $\mathcal{Y}$. Furthermore, $m_{\mathcal{X}}(A)$ equals $m_{\mathcal{Y}}\left(\varpi\big(A \cap \mathcal{X}_0^{\mathsf{c}}\big)\right)$. 
\end{lemma}

\begin{proof}

It will be easier first to prove the second statement about measures. For any set $B$ in $\mathcal{Y}$, a straightforward verification would tell us that we have the following equality of sets: \[\varpi^{-1}(B) \cap \mathcal{X}_0^{\mathsf{c}} = \varphi(B \cap \mathcal{Y}^{\mathsf{c}}).\] 
As a result, if $B$ is a measurable set in $\mathcal{Y}$, one can use the fact that $\mathcal{X}_0$ is a zero measure set to verify that
\begin{align}\label{eq:portmeasure}
m_{\mathcal{Y}}(B)= m_{\mathcal{X}}(\varpi^{-1}(B)) = m_{\mathcal{X}}(\varpi^{-1}(B) \cap \mathcal{X}_0^{\mathsf{c}}) = m_{\mathcal{X}}(\varphi(B \cap \mathcal{Y}_0^{\mathsf{c}})).
\end{align}
Since $\varpi$ maps $\mathcal{X}_0^{\mathsf{c}}$ to $\mathcal{Y}_0^{\mathsf{c}}$, note that $\varpi(A \cap \mathcal{X}_0^{\mathsf{c}})$ is a subset of $\mathcal{Y}_0^{\mathsf{c}}$. Setting $B=\varpi(A \cap \mathcal{X}_0^{\mathsf{c}})$ in equation (\ref{eq:portmeasure}) then tells us that  
\[
m_{\mathcal{Y}}\left(\varpi(A \cap \mathcal{X}_0^{\mathsf{c}})\right) = m_{\mathcal{X}}\left(A \cap \mathcal{X}_0^{\mathsf{c}}\right) = m_{\mathcal{X}}(A).
\]
To show that $\varpi\big(A \cap \mathcal{X}_0^{\mathsf{c}}\big)$ is a $P$-continuity set, we need to show that the measure of the boundary $\partial\left(\varpi\big(A \cap \mathcal{X}_0^{\mathsf{c}}\big)\right)$, with respect to $m_{\mathcal{Y}}$, equals zero. Setting $B$ equals $\partial\left(\varpi\big(A \cap \mathcal{X}_0^{\mathsf{c}}\big)\right)$ in equation (\ref{eq:portmeasure}), it is suffices to show that $\varphi\bigg(\partial\left(\varpi\big(A \cap \mathcal{X}_0^{\mathsf{c}}\big)\right) \cap \mathcal{Y}_0^{\mathsf{c}}\bigg)$ has zero measure in $\mathcal{X}$. We will show that 
\begin{align}\label{eq:boundary2}
\varphi\bigg(\partial\left(\varpi\big(A \cap \mathcal{X}_0^{\mathsf{c}}\big)\right) \cap \mathcal{Y}_0^{\mathsf{c}}\bigg) \stackrel{?}{\subset} \partial(A).
\end{align}  
This is sufficient for our purposes since $A$ is a $P$-continuity set in $\mathcal{X}$ and hence the measure $m_{\mathcal{X}}(\partial(A))$ equals zero.  

To show the inclusion in equation (\ref{eq:boundary2}), consider a point $x$ in the LHS of equation (\ref{eq:boundary2}). This point $x$ is of the form $\varphi(y)$, for some element $y$ in $\partial\left(\varpi\big(A \cap \mathcal{X}_0^{\mathsf{c}}\big)\right) \cap \mathcal{Y}_0^{\mathsf{c}}$. Since $y$ is in the boundary of $\varpi\big(A \cap \mathcal{X}_0^{\mathsf{c}}\big)$, there exists a sequence $\left\{y_n\right\}$ inside $\varpi\big(A \cap \mathcal{X}_0^{\mathsf{c}}\big)$ such that $y=\lim y_n$. If there existed any infinite subsequence of this sequence $\left\{y_n\right\}$ that lied entirely in $\mathcal{Y}_0$, then the limit point of that subsequence must also lie in $\mathcal{Y}_0$ since the set $\mathcal{Y}_0$ is closed. By assumption, however, $y$ belongs to $\mathcal{Y}_0^{\mathsf{c}}$. As a result, we may assume without loss of generality, that the $y_n$'s all belong to $\mathcal{Y}_0^{\mathsf{c}}$. Let $x_n$ denote $\varphi(y_n)$.  Since $x$ is in the image of $\varphi$, it must belong to $\mathcal{X}_0^{\mathsf{c}}$. Since $\mathcal{X}_0^{\mathsf{c}}$ is an open subset of $\mathcal{X}$, it is straightforward to check that $x = \lim x_n$. Since the $y_n$'s belong to $\varpi\big(A \cap \mathcal{X}_0^{\mathsf{c}}\big)$, the $x_n$'s belong to $A \cap \mathcal{X}_0^{\mathsf{c}}$, and hence belong to $A$. This shows that $x$ is a limit point in $A$.  To establish the inclusion in equation (\ref{eq:boundary2}), we need to show that $x$ is not an interior point of $A$. Suppose, for the sake of contradiction, that $x$ was an interior point of $A$. There must be an open set $U_x$ of $\mathcal{X}$ such that $x \in U_x \subset A$. Since $x$ belongs to $\mathcal{X}_0^{\mathsf{c}}$ and since $\mathcal{X}_0^{\mathsf{c}}$ is an open set in $\mathcal{X}$, the set $U_x \cap \mathcal{X}_0^{\mathsf{c}}$ is also an open set in $\mathcal{X}_0^{\mathsf{c}}$ containing $x$. That is, $x \in U_x \cap \mathcal{X}_0^{\mathsf{c}} \subset A \cap \mathcal{X}_0^{\mathsf{c}}$. Applying $\varpi$ to these inclusions, we obtain
\begin{align*}
y = \varpi(\varphi(y)) =  \varpi(x) \in \varpi\left(U_x \cap \mathcal{X}_0^{\mathsf{c}}\right) \subset \varpi\left(A \cap \mathcal{X}_0^{\mathsf{c}}\right) \subset \mathcal{Y}_0^{\mathsf{c}}.
\end{align*}

Since $\mathcal{Y}_0^{\mathsf{c}}$ is an open set of $\mathcal{Y}$ and since the restriction of $\varpi$ to $\mathcal{X}_0^{\mathsf{c}}$ is a homeomorphism, we must have that $\varpi\left(U_x \cap \mathcal{X}_0^{\mathsf{c}}\right)$ is an open set of $\mathcal{Y}$ containing $y$ and contained in $\varpi\left(A \cap \mathcal{X}_0^{\mathsf{c}}\right)$.  This observation lets us conclude that $y$ must be an interior point of $\varpi\left(A \cap \mathcal{X}_0^{\mathsf{c}}\right)$. This contradicts the fact that $y$ is in the boundary of $\varpi(A \cap \mathcal{X}_0^{\mathsf{c}})$, and hence not an interior point of $\varpi(A \cap \mathcal{X}_0^{\mathsf{c}})$. This establishes the inclusion in equation (\ref{eq:boundary2}) and thus proves that $\varpi(A)$ is a $P$-continuity set. The lemma follows. 
\end{proof}

\begin{lemma}\label{lem:transferequid}
Let $x$ be an element of $\mathcal{X}$. In addition to the conditions (\ref{item:gen1}), (\ref{item:gen2}) and (\ref{item:gen3}), suppose that the following conditions hold:
\begin{enumerate}[label=(\roman*), ref=\roman*]
\setcounter{enumi}{3}
\item\label{item:gen4} The set $\{\mathcal{T}_{\mathcal{Y}}^{[n]}\left(\varpi(x)\right)\}$ is equidistributed in $\mathcal{Y}$.
\item\label{item:gen5} There exists an $m$ such that the element $\mathcal{T}_{\mathcal{X}}^{[m]}(x)$ lies in $\mathcal{X}_0^{\mathsf{c}}$.
\end{enumerate}
Then, the sequence $\{\mathcal{T}_{\mathcal{X}}^{[n]}(x)\}$  is equidistributed in $\mathcal{X}$.
\end{lemma}

\begin{proof}
Since $\mathcal{T}_\mathcal{X}$ preserves $\mathcal{X}_0^{\mathsf{c}}$, the condition (\ref{item:gen5}) implies that for all  $n \geq m$, the element $\mathcal{T}_{\mathcal{X}}^{[n]}(x)$ lies in $\mathcal{X}_0^{\mathsf{c}}$. Since  equidistribution for a sequence is an asymptotic property, without loss of generality, we may assume that $m=0$ and that for all $n \geq 0$, the elements $\mathcal{T}_X^{[n]}(x)$ lie inside  $\mathcal{X}_0^{\mathsf{c}}$. In particular, we may assume without loss of generality that the point $x$ itself lies in the set $\mathcal{X}_0^{\mathsf{c}}$. By the Portmanteau Theorem (\cite[Theorem 2.1]{MR1700749}), to prove that the sequence $\{\mathcal{T}_{\mathcal{X}}^{[n]}(x)\}$  is equidistributed in $\mathcal{X}$, it is enough to show that for all $P$-continuity sets $A$ in $\mathcal{X}$, we have 
\begin{align}\label{eq:lemmaport}
	\lim_{M \rightarrow \infty} \dfrac{1}{M} \sum_{n=1}^M \mathds{1}_A\left(\mathcal{T}_{\mathcal{X}}^{[n]}(x)\right) \stackrel{?}{=} m_{\mathcal{X}}(A).
\end{align}
Here, $\mathds{1}_A$ is the indicator function for the set $A$. By our assumption, each of the points $\mathcal{T}_{\mathcal{X}}^{[n]}(x)$ belongs to the set $\mathcal{X}_0^{\mathsf{c}}$. As a result, for each term in the summand on the LHS of equation (\ref{eq:lemmaport}), we have 
\begin{align}\label{eq:transferequid1}
\mathds{1}_A\left(\mathcal{T}_{\mathcal{X}}^{[n]}(x)\right) = \mathds{1}_{A\cap\mathcal{X}_0^{\mathsf{c}}}\left(\mathcal{T}_{\mathcal{X}}^{[n]}(x)\right)	. 
\end{align}

Note that the restriction of $\varpi$ to $\mathcal{X}_0^{\mathsf{c}}$ is a homeomorphism, and in particular, a bijection onto its image. As a result, an element $z$ in $\mathcal{X}_0^{\mathsf{c}}$ belongs to $A\cap\mathcal{X}_0^{\mathsf{c}}$ if and only if its image $\varpi(z)$ belongs to $\varpi(A\cap\mathcal{X}_0^{\mathsf{c}})$. In other words, \[\mathds{1}_{A\cap\mathcal{X}_0^{\mathsf{c}}}(z) = \mathds{1}_{\varpi(A\cap\mathcal{X}_0^{\mathsf{c}})}(\varpi(z)).\] 
Applying this observation to the element $\mathcal{T}_{\mathcal{X}}^{[n]}(x)$  and using the commutative diagram in equation (\ref{commdiag:generalsetup}), we get 
\begin{align}\label{eq:transferequid2}
\mathds{1}_{A\cap\mathcal{X}_0^{\mathsf{c}}}\left(\mathcal{T}_{\mathcal{X}}^{[n]}(x)\right)	= \mathds{1}_{\varpi(A\cap\mathcal{X}_0^{\mathsf{c}})}\left(\mathcal{T}_{\mathcal{Y}}^{[n]}(\varpi(x))\right). 	
\end{align}
Using equations (\ref{eq:transferequid1}) and (\ref{eq:transferequid2}) along with Lemma \ref{lem:transferequid}, establishing the equation (\ref{eq:lemmaport}) reduces to establishing
\begin{align}\label{eq:porttoshow}
	\lim_{M \rightarrow \infty} \dfrac{1}{M} \sum_{n=1}^M \mathds{1}_{\varpi(A \cap \mathcal{X}_0^{\mathsf{c}})}\left(\mathcal{T}_{\mathcal{Y}}^{[n]}(\varpi(x))\right) \stackrel{?}{=} m_{\mathcal{Y}}\left(\varpi(A \cap \mathcal{X}_0^{\mathsf{c}})\right)
\end{align}
By Lemma \ref{lem:transferequid}, the set $\varpi(A \cap \mathcal{X}_0^{\mathsf{c}})$ is a $P$-continuity set. Equation (\ref{eq:porttoshow}) now follows by an application of the Portmanteau Theorem (\cite[Theorem 2.1]{MR1700749}) since the set $\{\mathcal{T}_{\mathcal{Y}}^{[n]}\left(\varpi(x)\right)\}$ is equidistributed in $\mathcal{Y}$. The lemma follows. 
\end{proof}

To complete the proof of the implication (\ref{it:prop7}) $\implies$ (\ref{it:prop2}), we will specialize the general setup to the situation we are interested in and verify the conditions (\ref{item:gen1}),(\ref{item:gen2}),(\ref{item:gen3a}),(\ref{item:gen3b}),(\ref{item:gen3c}),(\ref{item:gen4}) and (\ref{item:gen5}). We set
\begin{align*}
	& \mathcal{X} \coloneqq \left(\Z_p \times [0,1]\right)^r, \quad \mathcal{X}_0 \coloneqq \left(\Z_p\times \left\{0,1\right\}\right)^r, \quad \mathcal{T}_{\mathcal{X}} = T_3^r,  \quad m_X =\nu_3^r. \\ & \mathcal{Y} \coloneqq \left(\dfrac{\Z_p \times \mathbb{R}}{\mathbb{Z}}\right)^r, \quad \mathcal{Y}_0 \coloneqq \left(\dfrac{\Z_p \times \mathbb{Z}}{\mathbb{Z}}\right)^r, \quad \mathcal{T}_\mathcal{Y} = T_2^r, \quad m_{\mathcal{Y}} = \nu_2^r.
\end{align*}

The map $\varpi$ is the natural map $X_3^r \rightarrow X_2^r \cong X_1^r$ described in the introduction. Conditions (\ref{item:gen1}),(\ref{item:gen2}),(\ref{item:gen3a}) can be directly checked. Condition (\ref{item:gen3b}) follows from \cite[Appendix to Chapter 1, Proposition 1]{MR1760253}. Condition (\ref{item:gen3c}) follows from the observation that if there exists an integer $0\leq t_0 \leq p-1$ and a real number $0 < x <1$, then we have $0 < \dfrac{x+t_0}{p}< 1$.  Condition (\ref{item:gen4}) follows from the hypotheses of (\ref{it:prop4}). 

By the hypothesis in (\ref{it:prop4}), the set $\left\{T_2^{[n]}\big(\left(\gamma_1,x_1\right),\left(\gamma_2,x_2\right),\cdots,\left(\gamma_r,x_r\right)\big)\right\}$ is equidistributed in $X_2^r$. Under the continuous surjection $X_2^r \rightarrow (\mathbb{R}/\mathbb{Z})^r$, the image of this set must also be equidistributed in $(\mathbb{R}/\mathbb{Z})^r$ with respect to the Borel measure. If a sequence $\{z_n\}$ is equidistributed in $(\mathbb{R}/\mathbb{Z})^r$, then there must be infinitely many $n$ such that $z_n \notin \mathbb{Z}$. Using these observations, one can conclude that there must be at least once  $m$ such that the element $T_2^{[m]}\big(\left(\gamma_1,x_1\right),\left(\gamma_2,x_2\right),\cdots,\left(\gamma_r,x_r\right)\big)$ does not belong to $\mathcal{X}_0 \coloneqq \left(\Z_p\times \left\{0,1\right\}\right)^r$. This establishes condition (\ref{item:gen5}) and completes the proof of the implication (\ref{it:prop4}) $\implies$ (\ref{it:prop1}).
 
 \begin{remark}
 The idea for proving the implication (\ref{it:prop1}) $\iff$ (\ref{it:prop2}) in Proposition \ref{prop:nodependenceonrealcoordinates}, when dealing with skew product maps and Haar measures, can be traced back to a result of Furstenberg. See \cite[Proof of theorem 4.21]{MR2723325}.	
 \end{remark}

\section{Reduction to the $r=1$ case: Proof of Proposition \ref{prop:reductionstep}} \label{sec:proofpropreduction}

The main result in this section is Proposition \ref{prop:heckemap1pequiv}.  By setting $(\gamma_1,\cdots,\gamma_r)$ to equal $(\alpha\beta_1, \cdots, \alpha\beta_r)$ and using the implication (\ref{it:prop3}) $\iff$ (\ref{it:prop6}) of Proposition \ref{prop:nodependenceonrealcoordinates}, we see that statement \ref{it:prophecke1p1} of Proposition \ref{prop:heckemap1pequiv} is equivalent to statement \ref{it:propreduction1} of Proposition \ref{prop:reductionstep}. Once again, by using the implication (\ref{it:prop3}) $\iff$ (\ref{it:prop6}) of Proposition \ref{prop:nodependenceonrealcoordinates}, we see that statement \ref{it:prophecke1p2} of Proposition \ref{prop:heckemap1pequiv} is equivalent to statement \ref{it:propreduction2} of Proposition \ref{prop:reductionstep}. The proof of Proposition \ref{prop:reductionstep} would follow from Proposition \ref{prop:heckemap1pequiv}.

\begin{proposition}\label{prop:heckemap1pequiv}
	The following statements are equivalent:
	\begin{enumerate}[label=(\roman*), ref=(\roman*)]
		\item\label{it:prophecke1p1} The element $\big(\left(\gamma_1,0\right),\left(\gamma_2,0\right),\cdots,\left(\gamma_r,0\right)\big)$ in $X_2^r$ is a generic point for $T_2^r$.
		\item\label{it:prophecke1p2} For every $(m_1,\cdots,m_r) \in \mathbb{Z}^r \setminus \{(0,\cdots,0)\}$, the element $\bigg(\sum m_i\gamma_i,0\bigg)
		$ in $X_2$ is a generic point for $T_2$.
	\end{enumerate}
\end{proposition}

\begin{proof}

For each $\vec{m} = (m_1,\cdots,m_r)$ in $\mathbb{Z}^r \setminus \{(0,\cdots,0)\}$, we consider the map $\varpi_{\vec{m}}: X_2^r \rightarrow X_2$ given by 
\begin{align}\label{eq:formulaRQpmap}
\big(\left(\gamma_1,x_1\right),\left(\gamma_2,x_2\right),\cdots,\left(\gamma_r,x_r\right)\big) \mapsto \sum_{i=1}^r m_i(\gamma_i,x_i). 
\end{align}
Such a map is furnished since the map $(\R \times \Q_p)^r \rightarrow (\R\times\Q_p)^r$ given by equation (\ref{eq:formulaRQpmap}) preserves the lattice $(\Z[1/p])^r$. Observe that the continuous map $\varpi_{\vec{m}}: X_2^r \rightarrow X_2$ must be surjective since $\vec{m}$ is a non-zero vector. The measures $\nu_2^r$ and $\nu_2$ are the unique probability Haar measures on the compact topological groups $X_2^r$ and $X_2$ respectively. Using the fact that there exists a unique probability Haar measure on a compact topological group and the fact that the map $\varpi_{\vec{m}}: X_2^r \rightarrow X_2$ is a continuous surjective map, it is straightforward to imitate the argument in \cite[Page 20, point (4)]{MR648108} to conclude that the pushforward measure $\varpi_{\vec{m}}^*(\nu_2^r)$ on $X_2$ must coincide with the Haar measure $\nu_2$. The implication \ref{it:prophecke1p1} $\implies$ \ref{it:prophecke1p2} follows, just as in Section \ref{subsec:proofimplication}, by applying \cite[Theorem 3.6.1]{MR2267655}. \\ 

	To prove the reverse implication  \ref{it:prophecke1p2} $\implies$ \ref{it:prophecke1p1}, observe that, by considering the real and imaginary parts of a complex-valued function,  a sequence $\{z_n\}$ is equidistributed in $X_2^r$ if and only if for all complex valued functions  $f:X_2^r \rightarrow \mathbb{C}$, we have 
	\begin{align}\label{eq:equiddefrealfns}
	\lim_{M \rightarrow \infty}\dfrac{1}{M}\sum_{n=1}^M 	f(x_n) = \int f \ d\nu_2^r. 
	\end{align}

Since $X_2 \cong X_0 \cong \varprojlim \mathbb{R}/p^m\Z$, we can also consider the natural projection map $i_t:X_2 \rightarrow \mathbb{R}/p^t\Z$ . Finally, we will also consider the group isomorphism $e_t: \mathbb{R}/p^t\Z\xrightarrow {\cong} S^1$ given by the exponential map. Here are some preliminary observations:
\begin{align}\label{eq:complexSWapplicable}
\overline{e_t \circ i_t \circ \varpi_{\vec{m}}} &= e_t \circ i_t \circ \varpi_{-\vec{m}}, \\ \notag e_{t+1} \circ i_{t+1} \circ \varpi_{p \vec{m}} &= e_t \circ i_t \circ \varpi_{m},  \\ \notag (e_t \circ i_t \circ \varpi_{m_1}) (e_t \circ i_t \circ \varpi_{m_2})&= e_t \circ i_t \circ \varpi_{m_1+m_2}.
\end{align}

Consider the following collection of complex-valued functions on $X_2^r$:
\begin{align}\label{eq:complexfunctions}
	\mathcal{S} \coloneqq	\left\{c_0 + \sum_{{\vec{m}}} c_{\vec{m},t} {\cdot}\varphi_{\vec{m},t} : X_2^r \rightarrow \mathbb{C}, \text{where } c_0, c_{\vec{m},t} \in \mathbb{C}\right\}, 
	\end{align}
where we consider finite linear combinations in the summation given above in equation (\ref{eq:complexfunctions}) and where the map $\varphi_{\vec{m},t}:X_2^r \rightarrow \mathbb{C}^\times$ is the one described below:
	\begin{align}
		\xymatrix{
		X_2^r \ar[d]_{\varpi_{\vec{m}}}\ar@{..>}[rrr]^{\varphi_{\vec{m},t}}& & &  \mathbb{C}^\times, \\
		X_2 \ar[r]^{i_t \quad}& \mathbb{R}/p^{t}\Z \ar[r]^{e_t}_{\cong}& S^1 \ \ar@{^{(}->}[r]&  \mathbb{C}^\times \ar[u]^{=}
 		}
	\end{align}

	Using the observations in (\ref{eq:complexSWapplicable}), it is straightforward to verify that the hypotheses of the complex-form of Stone-Weierstrass theorem \cite[Theorem B.7]{MR2723325} is applicable to the set $\mathcal{S}$. Hence, the functions in $\mathcal{S}$ are dense in the set of complex-valued functions on $X_2^r$. Therefore, to etablish \ref{it:prophecke1p1}, i.e.~to show that the sequence $\left\{(T_2^r)^{[n]}\big(\left(\gamma_1,0\right),\left(\gamma_2,0\right),\cdots,\left(\gamma_r,0\right)\big)\right\}$ is equidistributed in $X_2^r$, it is enough to verify equation (\ref{eq:equiddefrealfns}) for the functions $\varphi_{\vec{m},t}$. Consider the RHS of (\ref{eq:equiddefrealfns}). Since the pushforward measure $\varpi_{\vec{m}}^*(\nu_2^r)$ on $X_2$ coincides with the Haar measure $\nu_2$, applying \cite[Theorem 3.6.1]{MR2267655} gives us the following equality:	
	\begin{align}\label{eq:RHSfinalprop}
	\int \phi_{\vec{m},t} \ d \nu_2^r= \int e_t \circ i_t \circ \varpi_{\vec{m}} \ d\nu_2^r	= \int e_t \circ i_t \ d \nu_2. 
	\end{align}
Since the maps $T_2^r$ on $X_2^r$ and $T_2$ on $X_2$ is given by the multiplication by $1/p$ map, we immediately have the following commutative diagram:
\begin{align*}
	\xymatrix{
	X_2^r \ar[d]_{\varpi_{\vec{m}}} \ar[r]^{T_2^r}& X_2^r \ar[d]^{\varpi_{\vec{m}}} \\
	X_2 \ar[r]^{T_2} & X_2.}
\end{align*}
Consider each term appearing in the summation on the LHS of (\ref{eq:equiddefrealfns}). We have 
\begin{align}\label{eq:LHSfinalprop}
 (e_t \circ i_t \circ \varpi_{\vec{m}}) (T_2^r)^{[n]}\big(\left(\gamma_1,0\right),\left(\gamma_2,0\right),\cdots,\left(\gamma_r,0\right)\big) = (e_t \circ i_t) \left((T_2)^{[n]}\left(\sum_{i=1}^r m_i\gamma_i,0\right)\right).
\end{align}
The desired equation in equation (\ref{eq:equiddefrealfns}) for the functions $\varphi_{\vec{m},t}$ now follows by using equations (\ref{eq:RHSfinalprop}) and (\ref{eq:LHSfinalprop}) along with the hypotheses \ref{it:prophecke1p2}. This proves the reverse implication \ref{it:prophecke1p2} $\implies$ \ref{it:prophecke1p1} and completes the proof of the proposition. 
\end{proof}

\begin{remark}
The possibility of this reduction itself is also briefly indicated in the original work of Ferrero--Washington \cite{MR528968}. In the reduction step from a general $r$ to the $r=1$ case, it seems to us that one cannot completely avoid using ``test functions'' in some form. However, in our proof, the exponential function is ``only'' used as a group homomorphism to identify $\R/p^t\Z$ with $\mathrm{S}^1$ inside $\mathbb{C}^\times$. What we really need is the existence of these test functions that separate points and preserve the group structure. What we want to emphasize is that this approach allows us to adopt a more abstract perspective where we can avoid any explicit computations of bounds involving test functions. 	
\end{remark}

\begin{remark}
We would like to highlight the contrast between working with the symbolic space $X_4$ verus working with the arithmetic description of the $p$-adic solenoid (as in $X_2$, $X_1$ or $X_0$). While the ergodic nature of 	the situation is captured by the symbolic space, as we remarked earlier, the reduction step is group-theoretic (or representation-theoretic) which necessitates working with the group structure of the $p$-adic solenoid. 
\end{remark}

\section{Acknowledgements}
We are very grateful to Manjunath Krishnapur, Fran\c{c}ois Ledrappier and Ruxi Shi for many enlightening discussions and helpful comments. We thank the referee for their insighful comments and suggestions.

This project began at the Mathematisches Forschungsinstitut Oberwolfach when BP was a Leibniz fellow and invited JL for a week-long stay at MFO. We thank MFO for its hospitality and for providing us with an excellent avenue to begin our collaboration. We are also grateful to the Heilbronn Institute for Mathematical Research (HIMR) and the UKRI/EPSRC Additional Funding Programme for Mathematical Sciences for funding JL's proposal \textit{New Perspectives in Iwasawa theory: Analytic side} via the Heilbronn Small Grant Scheme. JL is funded by ERC-Advanced Grant 833802 Resonances.  During the initial stages of this project, BP was a postdoctoral fellow at the National Center for Theoretical Sciences, Taiwan. BP would like to thank the faculty and staff affiliated with NCTS, especially Jungkai Chen, Ming-Lun Hsieh, Yng-Ing Lee and Peggy Lee, for their constant support and encouragement. BP's research is partially supported by the Infosys Young Investigator Award from the Infosys Foundation Bangalore along with the SERB-MATRICS grant MTR/2022/000244 and DST FIST program 2021 [TPN - 700661].

\bibliographystyle{abbrv}
\bibliography{bernoulli_shift_bib}

\Addresses

\end{document}